\documentclass[12pt]{article}
\usepackage{amsmath,amssymb,amsthm}

\newcommand{\gs}[1]{\textbf{#1}}

\newcommand{\fr}[2]{\frac{#1}{#2}}
\newcommand{\dfr}[2]{\dfrac{#1}{#2}}
\newcommand{\cd}{\cdot}
\newcommand{\cds}{\cdots}
\newcommand{\dsum}{\displaystyle \sum}

\renewcommand{\l}{\left}
\renewcommand{\r}{\right}

\newcommand{\vsv}{\vspace{5mm}}
\newcommand{\vsb}{\vspace{2mm}}

\newcommand{\q}{\quad}
\newcommand{\qq}{\qquad}

\newcommand{\maru}[1]{{\ooalign{\hfil#1\/\hfil\crcr
\raise.167ex\hbox{\mathhexbox20D}}}}

\newcommand{\ruby}[2]{%
 \leavevmode
 \setbox0=\hbox{#1}%
 \setbox1=\hbox{\tiny #2}%
 \ifdim\wd0>\wd1 \dimen0=\wd0 \else \dimen0=\wd1 \fi
 \hbox{%
   \kanjiskip=0pt plus 2fil
   \xkanjiskip=0pt plus 2fil
   \vbox{%
     \hbox to \dimen0{%
       \tiny \hfil#2\hfil}%
     \nointerlineskip
     \hbox to \dimen0{\mathstrut\hfil#1\hfil}}}}

\newcommand{\la}{\langle}
\newcommand{\ra}{\rangle}

\newcommand{\abs}[1]{\lvert{#1}\rvert}

\newcommand{\Z}{\mathbb{Z}}
\newcommand{\C}{\mathbb{C}}

\newcommand{\N}{\mathbb{N}}

\newcommand{\tensor}{\otimes}
\newcommand{\res}{\mathrm{Res}}

\newcommand{\vir}{\mathrm{Vir}}
\newcommand{\aut}{\mathrm{Aut}}

\renewcommand{\hom}{\mathrm{Hom}}
\newcommand{\ch}{\mathrm{ch}}
\newcommand{\SL}{\mathrm{SL}}
\newcommand{\id}{\mathrm{id}}
\newcommand{\ind}{\mathrm{Ind}}

\newcommand{\hf}{\fr{1}{2}}
\newcommand{\Span}{\mathrm{Span}}

\newcommand{\w}{\omega}
\newcommand{\vacuum}{\mathrm{1\hspace{-3.2pt}l}}
\newcommand{\vac}{\vacuum}

\newcommand{\fusion}{\boxtimes}

\makeatletter
\@addtoreset{equation}{section}

\theoremstyle{plain}
\newtheorem{thm}{Theorem}[section]
\newtheorem{prop}[thm]{Proposition}
\newtheorem{lem}[thm]{Lemma}
\newtheorem{cor}[thm]{Corollary}
\newtheorem{introthm}{Theorem}
\newtheorem{introcor}{Corollary}

\theoremstyle{definition}
\newtheorem{df}[thm]{Definition}
\newtheorem{df-lem}[thm]{Definition-Lemma}

\newtheorem{hypo}{Hypothesis}

\theoremstyle{remark}
\newtheorem{rem}[thm]{Remark}

\author{Hiroshi Yamauchi \vsb\\
  {\small \sl  Graduate School of Mathematics,}\\
  {\small \sl  University of Tsukuba, Ibaraki 305-8571, Japan}\\  
  {\small e-mail: {\tt hirocci@math.tsukuba.ac.jp}}
}

\title{2A-orbifold construction and the baby-monster \\
vertex operator superalgebra}

\date{}

\newcommand{\pf}{\gs{Proof:}\q }
\newcommand{\M}{\mathbb{M}}
\newcommand{\B}{\mathbb{B}}
\newcommand{\supp}{\mathrm{Supp}}

\newcommand{\com}{\mathrm{Com}}
\newcommand{\VB}{V\!\! B}

\begin{document}

\baselineskip 6mm

\maketitle

\begin{abstract}
In this article we prove that the full automorphism group of 
the baby-monster vertex operator superalgebra constructed by 
H\"{o}hn is isomorphic to $2\times \B$, where $\B$ is the baby-monster 
sporadic finite simple group, and determine irreducible modules for 
the baby-monster vertex operator algebra.
Our result has many corollaries.
In particular, we can prove that the $\Z_2$-orbifold 
construction with respect to a 2A-involution of the Monster 
applied to the moonshine vertex operator algebra $V^\natural$ yields 
$V^\natural$ itself again.
\end{abstract}

\section{Introduction}

The famous moonshine vertex operator algebra $V^\natural$ constructed 
by Frenkel-Lepowsky-Muerman \cite{FLM} is the first example of the 
$\Z_2$-orbifold construction of a holomorphic vertex operator algebra (VOA).
Let us explain a $\Z_2$-orbifold construction briefly.
Let $V$ be a holomorphic vertex operator algebra and $\sigma$ an 
involutive automorphism on $V$.
Then the fixed point subalgebra $V^{\la \sigma\ra}$ is a simple 
vertex operator algebra.
It is shown in \cite{DLM1} that there is a unique irreducible 
$\sigma$-twisted $V$-module $M$ and we have a decomposition
$M=M^0\oplus M^1$ into a direct sum of irreducible 
$V^{\la \sigma \ra}$-modules such that $M^0$ has an integral top weight.
Then a $\Z_2$-orbifold construction with respect to $\sigma\in \aut (V)$ 
is to construct a $\Z_2$-graded extension $W=V^{\la \sigma\ra}\oplus M^0$
of the fixed point subalgebra $V^{\la \sigma\ra}$ which is expected to 
be a holomorphic vertex operator algebra.

In FLM's construction, we take $V$ to be the lattice vertex operator 
algebra $V_\Lambda$ associated to the Leech lattice $\Lambda$ and 
the involution $\sigma$ is a natural lifting $\theta \in \aut (V_\Lambda)$ 
of the $(-1)$-isometry on $\Lambda$.
Denote by $V_\Lambda=V_\Lambda^+\oplus V_\Lambda^-$ the eigenspace
decomposition such that $\theta$ acts on $V_\Lambda^\pm$ as $\pm 1$.
Let $V_\Lambda^T$ be the unique irreducible $\theta$-twisted 
$V_\Lambda$-module.
Then there is a decomposition $V_\Lambda^T=(V_\Lambda^T)^+\oplus 
(V_\Lambda^T)^-$ such that the top weight of $(V_\Lambda^T)^+$ is integral.
Then the moonshine vertex operator algebra is defined by 
$V^\natural:= V_\Lambda^+\oplus (V_\Lambda^T)^+$ and 
it is proved in \cite{FLM} that $V^\natural$ forms a $\Z_2$-graded 
extension of $V_\Lambda^+$.
It is also proved in \cite{FLM} that the full automorphism group 
of the moonshine vertex operator algebra is the Monster sporadic 
finite simple group $\M$ by using Griess' result \cite{G}.

In the Monster, there are two conjugacy classes of involutions, 
the 2A-conjugacy class and the 2B-conjugacy class (cf.\ \cite{ATLAS}).
One can explicitly see the action of a 2B-involution on $V^\natural$ 
by FLM's construction.
But it is not clear to see the action of a 2A-involution on $V^\natural$
before Miyamoto.
In \cite{M1}, Miyamoto opened a way to study the action of 2A-involutions
of the Monster on the moonshine VOA by using a sub VOA isomorphic 
to the unitary Virasoro VOA $L(1/2,0)$.
Let us recall the definition of Miyamoto involutions.
Let $V$ be a simple VOA and $e\in V_2$ be a vector such 
that $e$ generates a sub VOA isomorphic to $L(1/2,0)$.
Such a vector $e$ is called conformal vector with central charge 1/2.
Since $V$ as a $\vir (e)\simeq L(1/2,0)$-module is completely 
reducible, we have a decomposition
$$
  V=V_e(0)\oplus V_e(1/2)\oplus V_e(1/16),
$$
where $V_e(h)$ denotes a sum of all irreducible $\vir (e)$-submodules
isomorphic to $L(1/2,h)$, $h=0,1/2,1/16$.
Then one can define a linear isomorphism $\tau_e$ on $V$ by
$$
  \tau_e :=\ 1 \q \text{on}\ \  V_e(0)\oplus V_e(1/2),\q 
  -1 \q \text{on}\ \  V_e(1/16).
$$
Then it is proved in \cite{M1} that $\tau_e$ defines an involution
of a VOA $V$ if $V_e(1/16)\ne 0$.
This involution is often called the (first) Miyamoto involution.
If $V_e(1/16)=0$, then one can define another automorphism on $V$ by
$$
  \sigma_e :=\ 1 \q \text{on}\ \  V_e(0),\q 
  -1 \q \text{on}\ \  V_e(1/2). 
$$
This involution is also called the (second) Miyamoto involution.
It is shown in \cite{C} and \cite{M1} that in the moonshine VOA every 
Miyamoto involution $\tau_e$ defines a 2A-involution of the Monster  
and the correspondence between conformal vectors and 2A-involutions is
one-to-one.
Therefore, in the study of 2A-involutions, it is very important 
to study conformal vectors with central charge 1/2.
Along this idea, C.H. Lam, H. Yamada and the author obtained 
an interesting achievement on 2A-involutions of the Monster in 
\cite{LYY}.

The main purpose of this paper is to study the $\Z_2$-orbifold 
construction of $V^\natural$ with respect to the Miyamoto involution 
and to prove that the 2A-orbifold construction applied to $V^\natural$ 
yields $V^\natural$ itself again.
Since a 2A-involution of the Monster is uniquely determined by 
a conformal vector $e$ of $V^\natural$ with central charge 1/2, 
we first study the commutant subalgebra of $\vir (e)$.
For a simple VOA $V$ and a conformal vector $e$ of $V$ with central
charge 1/2, set the space of highest weight vectors by 
$T_e(h):=\{ v\in V \mid e_{(1)}v\}$ for $h=0,1/2,1/16$.
Then we have decompositions $V_e(h)=L(1/2,h)\tensor T_e(h)$ and 
the commutant subalgebra $T_e(0)$ acts on $T_e(h)$ for $h=0,1/2,1/16$.
Since $L(1/2,0)$ has a $\Z_2$-graded extension $L(1/2,0)\oplus 
L(1/2,1/2)$, we can introduce a vertex operator superalgebra (SVOA)
structure on $T_e(0)\oplus T_e(1/2)$ (Theorem \ref{SVOA}) and 
its $\Z_2$-twisted module structure on $T_e(1/16)$ (Theorem 
\ref{twisted piece}).
It is easy to see that the one point stabilizer $C_{\aut (V)}(e)$ of a 
conformal vector $e$ naturally acts on the space of highest weight 
vectors $T_e(h)$.
If we take $V=V^\natural$, then $C_{\aut (V^\natural)}(e)$ is 
isomorphic to the 2-fold central extension $\la \tau_e \ra\cd \B$ of
the baby-monster sporadic finite simple group $\B$.
Therefore, the SVOA $T_e^\natural(0)\oplus T_e^\natural (1/2)$, where
we have set $V_e^\natural(h)=L(1/2,h)\tensor T_e^\natural(h)$ for 
$h=0,1/2,1/16$, affords a natural action of $\B$.
Motivated by this fact, H\"{o}hn first studied this SVOA in \cite{Ho1} 
and he called it the baby-monster SVOA.
Following him, we write $\VB^0:=T_e^\natural(0)$, $\VB^1:=T_e^\natural(1/2)$
and $\VB:=T_e^\natural(0)\oplus T_e^\natural(1/2)$.
It is proved in \cite{Ho2} that the full automorphism group of the even 
part $\VB^0$ of $\VB$ is exactly isomorphic to the baby-monster $\B$.
In this paper, we give a quite different proof of $\aut (\VB^0)\simeq \B$
based on a theory of simple current extensions.

In my recent work \cite{Y1} \cite{Y2}, a theory of simple current 
extensions of vertex operator algebras was developed and many useful 
results were obtained.
Using this theory, we determine the automorphism group of the commutant 
subalgebra $T_e(0)$ as follows:

\begin{introthm}
  Let $V$ be a holomorphic VOA and $e\in V$ a conformal vector with 
  central charge 1/2.
  Suppose the following:
  \\
  (a) $V_e(h)\ne 0$ for $h=0,1/2,1/16$, 
  \\
  (b) $V_e(0)$ and $T_e(0)$ are rational $C_2$-cofinite VOAs of CFT-type, 
  \\
  (c) $V_e(1/16)$ is a simple current $V^{\la \tau_e\ra}$-module,
  \\
  (d) $T_e(1/2)$ is a simple current $T_e(0)$-module,
  \\
  (e) $C_{\aut (V)}(e)/\la \tau_e\ra$ is a simple group or an odd group.
  \\
  Then 
  \\
  (i)\ $\aut (T_e(0))=C_{\aut (V)}(e)/\la \tau_e\ra$.
  \\
  (ii)\ The irreducible $T_e(0)$-modules are given by $T_e(0)$, 
  $T_e(1/2)$ and $T_e(1/16)$.
  \\
  (iii)\ The $\tau_e$-orbifold construction applied to $V$ yields 
  $V$ itself again.  
\end{introthm}

The assumptions (c) and (d) in the theorem above seem to be rather 
restrictive.
However, we prove that all the assumptions above hold if $V$ is 
the moonshine VOA.
We also present a refinement of Miyamoto's 
reconstruction of the moonshine VOA \cite{M5}.
Our refinement enable us to prove not only that the baby-monster SVOA 
$\VB$ satisfies all the assumptions above but also that we can 
construct the baby-monster SVOA $\VB$ without reference to $V^\natural$.
The main theorem of this paper is

\begin{introthm}
  Let $\VB=\VB^0\oplus \VB^1$ the simple SVOA obtained from $V^\natural$.
  \\
  (1) $\aut (\VB^0)= \B$ and $\aut (\VB)=2\times \B$.
  \\
  (2) There are exactly three irreducible $\VB^0$-modules, 
  $\VB^0$, $\VB^1$ and $\VB_T:=T_e^\natural(1/16)$.
  \\
  (3) The fusion rules for $\VB^0$-modules are as follows:
  $$
    \VB^1\times \VB^1=\VB^0,\q \VB^1\times \VB_T=\VB_T,\q 
    \VB_T\times \VB_T=\VB^0+\VB^1.
  $$
\end{introthm}

This theorem has many corollaries:

\begin{introcor}
  The irreducible 2A-twisted $V^\natural$-module has a shape
  $$
    L(1/2,1/2)\tensor \VB^0\oplus L(1/2,0)\tensor \VB^1
    \oplus L(1/2,1/16)\tensor \VB_T.
  $$
\end{introcor}

\begin{introcor}
  For any conformal vector $e\in V^\natural$ with central charge 1/2,
  there is no $\rho \in \aut (V^\natural)$ such that 
  $\rho (V_e^\natural(h))=V_e^\natural(h)$ for $h=0,1/2,1/16$ and 
  $\rho|_{(V^\natural)^{\la \tau_e\ra}}=\sigma_e$.
\end{introcor}

\begin{introcor}
  The 2A-orbifold construction applied to the moonshine VOA $V^\natural$
  yields $V^\natural$ itself again.
\end{introcor}

At the end of this paper, we give characters of $\VB^0$-modules and their
modular transformation laws.
Surprisingly, we find that the fusion algebra and the modular transformation
laws for the baby-monster VOA is canonically isomorphic to those of the 
Ising model $L(1/2,0)$.

\section{Simple current extension}

Let $V$ be a simple vertex operator algebra (VOA).
We recall a definition of a fusion product of $V$-modules. 

\begin{df}
  Let $M^1$, $M^2$ be $V$-modules.
  A {\it fusion product} for the ordered pair $(M^1,M^2)$ is a pair 
  $(M^1\fusion_{V}M^2 ,F(\cd,z))$ consisting of a $V$-module 
  $M^1 \fusion_{V} M^2$ and a $V$-intertwining operator 
  $F(\cd,z)$ of type $M^1\times M^2\to M^1\fusion_{V} M^2$ 
  satisfying the following universal property:
  For any $V$-module $W$ and any $V$-intertwining operator 
  $I(\cd,z)$ of type $M^1\times M^2 \to W$, there exists a
  unique $V$-homomorphism $\psi$ from $M^1\fusion_{V} 
  M^2$ to $W$ such that $I(\cd,z)=\psi F(\cd,z)$.
  We usually denote the pair $(M^1\fusion_V M^2, F(\cd,z))$ 
  simply by $M^1\fusion_V M^2$.
\end{df}

A theory of fusion products has been greatly developed by 
Huang-Lepowsky \cite{HL1}-\cite{HL4} and Huang \cite{H1}-\cite{H4}
(see also \cite{DLM2} \cite{Li3}), and it is proved that 
if $V$ is rational then a fusion product of any two $V$-modules always
exists (cf.\ \cite{HL3} \cite{HL4} \cite{Li3}) and if $V$ is also 
$C_2$-cofinite and of CFT-type then the fusion product satisfies the 
associativity (cf.\ \cite{H1} \cite{H4} \cite{DLM2}).
Therefore, the theory of fusion products is a powerful tool to 
study a rational $C_2$-cofinite vertex operator algebra of CFT-type.
Among modules for such a vertex operator algebra, simple current 
modules have a special importance.

\begin{df}
  An irreducible $V$-module $U$ is called a {\it simple 
  current} if it satisfies: For any irreducible $V$-module 
  $W$, the fusion product $U\fusion_{V} W$ is also irreducible.
\end{df}

In this paper we mainly consider the following extensions of 
vertex operator algebras. 

\begin{df}
  Let $V^0$ be a simple rational VOA and $D$ a finite abelian group.
  Let $\{ V^\alpha \mid \alpha \in D\}$ be a set of inequivalent
  irreducible $V^0$-modules.
  A {\it $D$-graded extension $V_D$ of $V^0$} is a simple VOA 
  $V_D=\oplus_{\alpha \in D} V^\alpha$ which extends the VOA structure 
  of $V^0$ with the grading structure 
  $Y_{V_D}(x^\alpha,z) x^\beta \in V^{\alpha +\beta}((z))$ 
  for any $x^\alpha \in V^\alpha$, $x^\beta\in V^\beta$.
  A $D$-graded extension $V_D$ is called a {\it $D$-graded simple 
  current extension of $V^0$} if all $V^\alpha$, $\alpha \in D$, 
  are simple current $V^0$-modules.
  In the case of $D=\Z_2=\{ 0,1\}$ and $V_D=V^0\oplus V^1$ is a simple 
  vertex operator superalgebra with even part $V^0$ and odd part $V^1$, 
  we call $V_D$ a {\it simple current super-extension of $V^0$} 
  if $V^1$ is a simple current $V^0$-module.
\end{df}

\begin{rem}
  Let $D^*$ be the dual group of $D$.
  By definition, there is a natural action of $D^*$ on $V_D$
  defined as $\chi|_{V^\alpha}=\chi (\alpha)\cd \id_{V^\alpha}$ for 
  $\alpha \in D$ and $\chi \in D^*$.
\end{rem}

We have the following uniqueness property of a simple 
current extension.

\begin{lem}\label{uniqueness0}
  (\cite{DM2})
  (i) Let $V_D=\oplus_{\alpha \in D} V^\alpha$ be a $D$-graded extension
  of $V^0$.
  If the space of $V^0$-intertwining operators of type $V^\alpha \times 
  V^\beta\to V^{\alpha+\beta}$ is one-dimensional for all $\alpha,\beta
  \in D$, then the VOA structure on $V_D$ is unique over $\C$.
  In particular, if $V_D$ is a $D$-graded simple current extension
  of $V^0$, then its VOA structure is unique over $\C$.
  \\
  (ii) The SVOA structure on a simple current super-extension is unique over 
  $\C$.
\end{lem}

In general, it is a difficult problem to determine whether a given module 
is a simple current or not.
However, the following lemma provides us a simple criterion.

\begin{lem}\label{criterion}
  (\cite{Y2}) 
  Let $V$ be a simple rational $C_2$-cofinite VOA of CFT-type and $U$ a 
  $V$-module.
  If there is a $V$-module $W$ such that the fusion rule $U \fusion_V W =V$ 
  holds in the fusion algebra for $V$, then $U$ (and also $W$) is a simple 
  current $V$-module.
\end{lem}

\pf
By the assumption, we can use the results in \cite{H4} so that 
the fusion algebra for $V$ is a commutative associative algebra over $\N$ 
with the unit element $V$.
Let $X$ be an irreducible $V$-module.
Then $U \fusion_V X$ is also a $V$-module and is a direct sum of irreducible
$V$-modules.
Let $U \fusion_V X=\oplus_{i\in I} T^i$ be a decomposition into a 
direct sum of irreducible $V$-modules.
First, we show that $W \fusion_V T^i\ne 0$ for all $i\in I$.
Assume that there is an $i_0\in I$ such that $W \fusion_V T^{i_0} =0$.
Then by multiplying $U$ in the fusion algebra we obtain $0=U \fusion_V
(W \fusion_V T^{i_0})=(U \fusion_V W)\fusion_V T^{i_0}= V \fusion_V T^{i_0}
=T^{i_0}$, a contradiction.
Therefore, $W \fusion_V T^i\ne 0$ for all $i \in I$.
Then by multiplying $W$ to $U\fusion_V X$ we get $X= V\fusion_V X
=(W \fusion_V U) \fusion_V X= W\fusion_V (U\fusion_V X) 
= W \fusion_V \oplus_{i\in I} T^i = \oplus_{i\in I} W\fusion_V T^i$.
Therefore, the cardinality of the index set $I$ is 1 and hence 
$U\fusion_V X$ is an irreducible $V$-module.
Thus $U$ is a simple current $V$-module.
\qed

The representation theory of simple current extensions is studied 
in many papers (cf.\ \cite{DLM3} \cite{L} \cite{SY} \cite{Y1} \cite{Y2}) 
and it is shown in \cite{L} \cite{Y1} \cite{Y2} that every simple current 
extension of a simple rational $C_2$-cofinite VOA of CFT-type is also 
rational and $C_2$-cofinite.
We review some results from \cite{Y1} and \cite{Y2} which we will need later.

Let $V_D=\oplus_{\alpha \in D} V^\alpha$ be a $D$-graded simple current 
extension of a simple rational $C_2$-cofinite CFT-type VOA $V^0$.
Since the fusion algebra for $V^0$ is associative, we can adopt the following
definition.

\begin{df}
  Let $W$ be an irreducible $V^0$-module $W$.
  A subset 
  $D_W:=\{ \alpha \in D \mid V^\alpha \fusion_V W\simeq W\}$
  forms a subgroup of $D$.
  We call $D_W$ the {\it stabilizer} of $W$.
\end{df}

\begin{lem}\label{SY}
  (\cite{SY})
  Let $M$ be a $V_D$-module and $W$ an irreducible $V^0$-submodule of $M$.
  Then $V^\alpha\cd W$ is also a non-trivial irreducible $V^0$-submodule 
  of $M$, where $V^\alpha\cd W$ denotes a linear space spanned by 
  elements $a_{(n)} w=\res_z z^n Y_{M}(a,z)w$ with $a\in V^\alpha$, 
  $w\in W$ and $n\in \Z$. 
  The stabilizer $D_W$ is determined independently of a choice of an 
  irreducible $V^0$-submodule $W$ if $M$ is an irreducible $V_D$-module.
\end{lem}

By the lemma above, we introduce the following notion.

\begin{df}
  An irreducible $V_D$-module $M$ is called {\it $D$-stable} if
  $D_W=0$ for some irreducible $V^0$-submodule $W$ of $M$.
\end{df}

Among $V_D$-modules, $D$-stable modules enjoy nice properties.

\begin{prop}\label{uniqueness1}
  (\cite{SY} \cite{Y1})
  Let $M$ be an irreducible $D$-stable $V_D$-module.
  Then the $V_D$-module structure on a $V^0$-module $M$ is unique
  over $\C$.
\end{prop}

\begin{thm}\label{induced module}
  (\cite{Y1}, Induced modules) 
  Let $W$ be an irreducible $V^0$-module.
  Then there exists a unique $\chi \in D^*$ such that $W$ is contained 
  in an irreducible $\chi$-twisted $V_D$-module.
  If $D_W=0$, then it is given by the induced module
  $$
    \ind_{V^0}^{V_D} W:=\bigoplus_{\alpha \in D} V^\alpha \fusion_{V^0} W .
  $$
  Moreover, any irreducible $\chi$-twisted $V_D$-module containing $W$ 
  as a $V^0$-submodule is isomorphic to $\ind_{V^0}^{V_D} W$ above.
\end{thm}

For a later purpose, we give a detailed description of the theorem above
in the case of simple current super-extensions.

\begin{thm}\label{super induced module}
  (\cite{Y2}) 
  Let $V=V^0\oplus V^1$ be a simple current super-extension of a simple 
  rational $C_2$-cofinite VOA $V^0$ of CFT-type.
  For an irreducible $V^0$-module $W$, 
  \\
  (i) If $V^1\fusion_{V^0} W\not\simeq W$, then $W$ is uniquely lifted to be 
  an irreducible untwisted or $\Z_2$-twisted 
  $V$-module $W\oplus (V^1\fusion_{V^0} W)$.
  \\
  (ii) If $V^1\fusion_{V^0} W\simeq W$, then there exist exactly two 
  inequivalent irreducible $\Z_2$-twisted $V$-module structure on $W$.
  If we write one of them by $W^+$, then the other one is given as 
  the $\Z_2$-conjugate $V$-module of $W^+$.
\end{thm}

\begin{thm}\label{lifting}
  (\cite{SY} \cite{Y1}, Lifting property of intertwining operators)
  Let $M^i$, $i=1,2,3$, be irreducible $D$-stable $V_D$-modules 
  and let $W^i$ be irreducible $V^0$-submodules of $M^i$ for $i=1,2,3$.
  Then for a $V^0$-intertwining operator $I(\cd,z)$ of type 
  $W^1\times W^2\to W^3$ there is a $V_D$-intertwining operator 
  $\tilde{I}(\cd,z)$ of type $M^1\times M^2\to M^3$ such that 
  $\tilde{I}(\cd,z)|_{W^1\tensor W^2}=I(\cd,z)$.
  Therefore, there is a natural linear isomorphism
  $$
    \binom{M^3}{M^1\ M^2}_{V_D} \simeq \bigoplus_{\alpha \in D} 
    \binom{\ V^\alpha \fusion_{V^0} W^3\ }{W^1\qq W^2}_{V^0},
  $$
  where $\binom{M^3}{M^1\ M^2}_{V_D}$ denotes the space of 
  $V_D$-intertwining operators of type $M^1\times M^2\to M^3$.
\end{thm}

\begin{cor}\label{cor:2.13}
  Let $E$ be a subgroup of $D$ and fix a coset decomposition
  $D=\sqcup_{i=1}^r (\alpha_i+E)$ of $D$.
  Then $V_E:=\oplus_{\alpha \in E} V^\alpha$ is an $E$-graded 
  simple current extension of $V^0$ and $V_{E+\alpha_i}:= \oplus_{\beta\in E}
  V^{\alpha_i+\beta}$ are simple current $V_E$-modules.
  Therefore, we can view $V_D=\oplus_{i=1}^r V_{E+\alpha_i}$ as a $D/E$-graded 
  simple current extension of $V_E$.
\end{cor}

Let us consider automorphisms on $V_D$.
Let $\sigma \in \aut (V^0)$ and $(X,Y_X(\cd,z))$ a $V^0$-module.
We can define the {\it $\sigma$-conjugate module} $X^\sigma$ as follows.
As a vector space, we set $X^\sigma= X$ and the vertex operator map 
on $X^\sigma$ is defined by $Y_{X^\sigma}(a,z):=Y_X(\sigma a,z)$ for
$a\in V^0$.
It is clear that $X^\sigma$ is irreducible if and only if $X$ is 
an irreducible $V^0$-module.
An irreducible $V^0$-module $X$ is called {\it $\sigma$-stable} if 
its $\sigma$-conjugate $X^\sigma$ is isomorphic to $X$ as a $V^0$-module.

\begin{rem}
  The definition of $D$-stable $V_D$-modules comes from the fact 
  that if $M$ is an irreducible $D$-stable $V_D$-module then
  its conjugate $M^\chi$ is isomorphic to $M$ as a $V_D$-module 
  for any $\chi \in D^*$.
\end{rem}

We can construct the {\it $\sigma$-conjugate} $(V_D)^\sigma$ 
of $V_D$ as follows.
Let $Y_{V_D}(\cd,z)$ be the vertex operator map on $V_D$.
By definition, there are canonical linear isomorphisms $\psi_\alpha:
V^\alpha \to (V^\alpha)^\sigma$, $\alpha \in D$, such that
$$
  Y_{(V^\alpha)^\sigma}(x^0,z)\psi_\alpha = \psi_\alpha 
  Y_{V^\alpha}(\sigma x^0,z)
$$
for all $x^0\in V^0$.
Then we define the vertex operator map $Y_{V_D}^\sigma(\cd,z)$ 
on $(V_D)^\sigma=\oplus_{\alpha \in D}(V^\alpha)^\sigma$ by
$$
  Y_{V_D}^\sigma(\psi_\alpha x^\alpha,z) \psi_\alpha x^\beta
  := \psi_{\alpha+\beta} Y_{V_D}(x^\alpha,z) x^\beta
$$
for $x^\alpha \in V^\alpha$, $x^\beta \in V^\beta$.
Then one can easily verify that $((V_D)^\sigma, Y_{V_D}^\sigma(\cd,z))$
also forms a $D$-graded simple current extension of $V^0$.
The following lifting property is established in \cite{Sh} by using 
the uniqueness of VOA structure on $V_D$.

\begin{thm}\label{Sh}
  (\cite{Sh}, Lifting property of automorphisms)
  Let $\sigma \in \aut (V^0)$ such that $(V_D)^\sigma \simeq V_D$ 
  as a $V^0$-module.
  Then there is a lifting $\tilde{\sigma}\in \aut (V_D)$ such that
  $\tilde{\sigma} V^0=V^0$ and $\tilde{\sigma}|_{V^0}=\sigma$.
  The lifting $\tilde{\sigma}$ is unique up to multiples of 
  elements in $D^*\subset \aut (V_D)$.
  This assertion still holds if $D=\Z_2$ and $V_D=V^0\oplus V^1$ is a 
  simple current super-extension of $V^0$.
\end{thm}

As we have seen above, simple current extensions have many good 
properties.
At last of this section, we present the following extension property
of simple current extensions.

\begin{thm}\label{extension}
  (\cite{Y2}, Extension property of simple current extensions) 
  Let $V^{(0,0)}$ be a simple rational $C_2$-cofinite VOA of CFT-type, 
  and let $D_1$, $D_2$ be finite abelian groups.
  Assume that we have a set of inequivalent irreducible simple current 
  $V^{(0,0)}$-modules $\{ V^{(\alpha,\beta)} \mid (\alpha,\beta) \in 
  D_1\oplus D_2\}$ with $D_1\oplus D_2$-graded  fusion rules 
  $V^{(\alpha_1,\beta_1)} \fusion_{V^{(0,0)}} V^{(\alpha_2,\beta_2)}=  
  V^{(\alpha_1 +\alpha_2,\beta_1+\beta_2)}$ for any 
  $(\alpha_1,\beta_1), (\alpha_2,\beta_2) \in D_1\oplus D_2$.
  \\
  (i) Further assume that all $V^{(\alpha,\beta)}$, $(\alpha,\beta) \in D_1\oplus 
  D_2$, have integral top weights and we have $D_1$- and $D_2$-graded simple 
  current extensions $V_{D_1}=\oplus_{\alpha \in D_1} V^{(\alpha,0)}$ and 
  $V_{D_2}= \oplus_{\beta \in D_2} V^{(0,\beta)}$.
  Then $V_{D_1\oplus D_2}:=\oplus_{(\alpha,\beta) \in D_1\oplus D_2} 
  V^{(\alpha,\beta)}$ possesses a unique structure of a simple vertex operator 
  algebra as a $D_1\oplus D_2$-graded simple current extension of $V^{(0,0)}$.
  \\
  (ii) In the case of $D_2=\Z_2=\{ 0,1\}$, further assume that 
  all $V^{(\alpha,0)}$, $\alpha \in D_1$, have integral top weight, 
  all $V^{(\alpha,1)}$, $\alpha \in D_1$, have  half-integral top weight, 
  $V_{D_1}=\oplus_{\alpha \in D_1} V^{(\alpha,0)}$ is a $D_1$-graded simple 
  current extension of $V^{(0,0)}$, and $V_{D_2}=V^{(0,0)}\oplus V^{(0,1)}$ 
  is a simple current super-extension of $V^{(0,0)}$. 
  Then $V_{D_1\oplus D_2}=\oplus_{(\alpha,\beta)\in D_1\oplus D_2} 
  V^{(\alpha,\beta)}$ has a unique structure of a simple vertex operator 
  superalgebra with even part $\oplus_{\alpha \in D_1} V^{(\alpha,0)}$ and 
  odd part $\oplus_{\beta \in D_1}V^{(\beta,1)}$ as a simple current 
  super-extension of $V_{D_1}$.
\end{thm}

\section{Miyamoto involution and its centralizer}

Let us denote by $L(c,h)$ the irreducible highest weight module for the 
Virasoro algebra with central charge $c$ and highest weight $h$.
It is shown in \cite{FZ} that $L(c,0)$ has a structure of a 
simple VOA.
Here we consider the first unitary Virasoro VOA $L(1/2,0)$.
It is proved in \cite{DMZ} \cite{W} that $L(1/2,0)$ is a rational 
$C_2$-cofinite VOA of CFT-type and has exactly three
irreducible modules, $L(1/2,0)$, $L(1/2,1/2)$ and $L(1/2,1/16)$.
Their fusion rules have also been computed and are as follows:

\begin{equation}\label{fusion ising}
\begin{array}{l}
  L(1/2,1/2)\times L(1/2,1/2)=L(1/2,0),\q 
  \vsb\\
  L(1/2,1/2)\times L(1/2,1/16)=L(1/2,1/16),
  \vsb\\
  L(1/2,1/16)\times L(1/2,1/16)=L(1/2,0)+L(1/2,1/2).
\end{array}
\end{equation}

First, we present an explicit realization of $L(1/2,0)$ and its modules.

\subsection{Ising model}\label{Ising}

In this section we give an explicit construction of the Ising
model SVOA $L(1/2,0)\oplus L(1/2,1/2)$ and its $\Z_2$-twisted modules
$L(1/2,1/16)^\pm$.
This construction is well-known and the most of contents in this section
can be found in \cite{FFR}, \cite{FRW} and \cite{KR}. 

Let $\mathcal{A}_\psi$ be the algebra generated by $\{ \psi_k \mid k\in 
\Z+\fr{1}{2}\}$ subject to the defining relation
$$
\begin{array}{ll}
  [\psi_m,\psi_n]_+:=\psi_m\psi_n+\psi_n\psi_m =\delta_{m+n,0},
  & m,n \in \Z +\fr{1}{2},
\end{array}
$$
and denote a subalgebra of $\mathcal{A}_\psi$ generated by $\{ 
\psi_k \mid k\in \Z +\fr{1}{2}, k>0\}$ by $\mathcal{A}_\psi^+$.
Let $\C \vac$ be a trivial $\mathcal{A}_\psi^+$-module.
Define a canonical induced $\mathcal{A}_\psi$-module $M$ by
$$
  M:=\mathrm{Ind}_{\mathcal{A}_\psi^+}^{\mathcal{A}_\psi} \C \vac
  = \mathcal{A}_\psi \tensor_{\mathcal{A}_\psi^+} \C \vac .
$$
Consider the generating series
$$
  \psi (z):=\dsum_{n\in\Z}\psi_{n+\fr{1}{2}}z^{-n-1}.
$$
Since $[\psi (z),\psi (w)]_+ = z^{-1}\delta (\fr{w}{z})$, $\psi (z)$ 
is local with itself.
So we can consider a subalgebra of a local system on 
$M$ generated by $\psi (z)$ and $I(z)=\id_M$.
By a direct calculation, one sees that the component operators of 
the generating series 
$$
  \fr{1}{2}\psi (z)\circ_{-2}\psi (z)
  :=\fr{1}{2}\res_{z_0}\{ (z_0-z)^{-2}\psi (z_0)\psi(z)+(-z+z_0)^{-2}
  \psi (z_0) \psi (z)\}
$$
defines a representation of the Virasoro algebra on $M$ with central 
charge 1/2, where $\circ_n$ denotes the $n$-th normal ordered product defined 
in \cite{Li1}.
It follows from \cite{KR} that $M$ as a $\vir$-module is isomorphic to
$L(1/2,0)\oplus L(1/2,1/2)$.
Therefore, there is a unique simple vertex operator superalgebra 
structure on $M$ such that $Y_M(\psi_{-\fr{1}{2}}\vac,z)=\psi (z)$.

\begin{thm} 
  On $M$, there is a unique simple vertex operator superalgebra structure
  $(M, Y_M(\cd,z), \vac,\fr{1}{2}\psi_{-\fr{3}{2}} \psi_{-\fr{1}{2}}\vac )$ 
  such that $Y_M(\psi_{-\fr{1}{2}}\vac,z)=\psi (z)$.
\end{thm}

Another unitary  $\vir$-module $L(1/2,1/16)$ is realized as follows.
Let $\mathcal{A}_\phi$ be the algebra generated by $\{ \phi_n \mid
 n\in\Z\}$ with defining relation
$$
  [\phi_m,\phi_n]_+=\delta_{m+n,0},\q m,n\in\Z .
$$
Let $\mathcal{A}_\phi^+$ be a subalgebra of $\mathcal{A}_\phi$
generated by $\{ \phi_n | n>0\}$ and let $\C v_0$ be a trivial 
one-dimensional $\mathcal{A_\phi^+}$-module. 
Then set $N=\mathrm{Ind}_{\mathcal{A}_\phi^+}^{\mathcal{A}_\phi}\C
v_0$ as we did previously. 
We can find an action of the Virasoro algebra on $N$. 
Consider the generating series
$$
  \phi (z):= \dsum_{n\in\Z}\phi_n z^{-n-\fr{1}{2}}.
$$
By direct calculations one can show that $\phi (z)$ is local with
itself.
Consider a local system on $N$ containing $\phi (z)$.
Since the powers of $z$ in $\phi (z)$ lie in $\Z +\fr{1}{2}$, we have
to use the $\Z_2$-twisted normal ordered product in \cite{Li2}. 
Define a generating series $L(z)$ of operators
on $N$ by
$$
\begin{array}{ll}
     L(z)
     &:= \phi (z)\circ_{-2} \phi(z)
     \vsb\\
     &= \dfr{1}{2}\res_{z_0} \res_{z_1} z_0^{-2} \l( 
       \dfr{z_1-z_0}{z} \r)^{\fr{1}{2}} 
     \vsb\\
     &\qq \times 
     \l\{ z_0^{-1}\delta\l(\dfr{z_1-z}{z_0}\r) \phi (z_1)  
       \phi (z) + z_0^{-1}\delta\l(\dfr{-z+z_1}{z_0}\r) \phi (z) 
       \phi (z_1)\r\} ,
     \vsb
\end{array}
$$
where $\circ_{n}$ above denotes the $n$-th normal ordered product 
in a $\Z_2$-twisted local system on $N$ (cf.\ \cite{Li2}).
Then by a direct computation we find that the component operators of 
$L(z)$ defines a representation of the Virasoro algebra on $N$ with 
central charge 1/2.
Set $v_{1/16}^\pm:=\phi_0\vac \pm (1/\sqrt{2})\vac$.
Then $v_{1/16}^\pm$ are highest weight vectors for the Virasoro 
algebra and the following decomposition is shown in \cite{KR}:
$$
  N= L(1/2,1/16)^+\oplus L(1/2,1/16)^-,
$$
where $L(1/2,1/16)^\pm$ are highest weight $\vir$-module generated by 
$v_{1/16}^\pm$, respectively.

\begin{thm}\label{1/16pm}
  The following $\Z_2$-twisted Jacobi identity holds on $N$:
  $$
  \begin{array}{l}
    z_0^{-1}\delta \l(\dfr{z_1-z_2}{z_0}\r) Y_N(a,z_1) Y_N(b,z_2)
    -(-1)^{\varepsilon (a,b)} z_0^{-1} \l(\dfr{-z_2+z_1}{z_0}\r) 
      Y_N(b,z_2) Y_N(a,z_1) 
    \vsb\\
    =z_1^{-1}\l(\dfr{z_2+z_0}{z_1}\r) \l(\dfr{z_2+z_0}{z_1}
      \r)^{\varepsilon(a,a)/2} Y_N(Y_M(a,z_0)b,z_2),
  \end{array}
  $$
  where $a,b\in M=L(1/2,0)\oplus L(1/2,1/2)$ and $\varepsilon(\cd,\cd)$ 
  denotes the standard parity function.
  Therefore, the vertex operator map $Y_N(\cd,z)$ defines inequivalent 
  irreducible $\Z_2$-twisted $L(1/2,0)\oplus L(1/2,1/2)$-module structures
  on $L(1/2,1/16)^\pm$.
\end{thm}

\subsection{Miyamoto involution}

Let $(V,Y_V(\cd,z),\vac,\w)$ be a VOA.
A vector $e\in V$ is called a {\it conformal vector} if its component operators 
$Y_V(e,z)=\sum_{n\in \Z} e_{(n)}z^{-n-1}=\sum_{n\in \Z} L^e(n)z^{-n-2}$ 
generate a representation of the Virasoro algebra on $V$:
$$
  [L^e(m), L^e(n)]=(m-n)L^e(m+n)+\delta_{m+n,0} \dfr{m^3-m}{12} c_e.
$$
The scalar $c_e$ is called {\it central charge} of a conformal vector $e$.
We denote by $\vir (e)$ the sub VOA generated by $e$.
If $\vir (e)$ is a rational VOA, then $e$ is called a {\it rational conformal 
vector}.
A decomposition $\w=e+(\w-e)$ is called {\it orthogonal} if both $e$ and 
$\w-e$ are conformal vectors and their component operators are mutually 
commutative.

Now assume that $e\in V$ is a rational conformal vector with central 
charge 1/2.
Then $\vir (e)$ is isomorphic to $L(1/2,0)$. 
Since $L(1/2,0)$ is rational, we can decompose $V$ into a direct sum of 
irreducible $\vir (e)$-modules as follows:
$$
  V=V_e(0)\oplus V_e(1/2) \oplus V_e(1/16),
$$
where $V_e(h)$, $h\in \{ 0,1/2,1/16\}$, denotes the sum of all 
irreducible $\vir (e)$-submodules of $V$ isomorphic to $L(1/2,h)$.
By the fusion rules \eqref{fusion ising}, we have the following grading 
structure (cf.\ \cite{M1}):
$$
\begin{array}{l}
  V_e(0)\cd V_e(h)\subset V_e(h),\q h=0,1/2,1/16,\q
  V_e(1/2)\cd V_e(1/2)\subset V_e(0),
  \vsb\\
  V_e(1/2)\cd V_e(1/16)\subset V_e(1/16),\q 
  V_e(1/16)\cd V_e(1/16)\subset V_e(0)\oplus V_e(1/2).
\end{array}
$$
Therefore, if $V_e(1/16)\ne 0$, then the linear map
$$
  \tau_e :=\ 1 \q \text{on}\ \  V_e(0)\oplus V_e(1/2),\q 
  -1 \q \text{on}\ \  V_e(1/16)
$$
defines an involutive automorphism on $V$ (cf.\ \cite{M1}).
We call $\tau_e$ the {\it first Miyamoto involution} or simply {\it 
Miyamoto involution} associated to a conformal vector $e$.
If $V_e(1/16)=0$, then we can also define another involution as follows 
(cf.\ \cite{M1}):
$$
  \sigma_e:=\ 1 \q \text{on}\ \ V_e(0),\q
  -1 \q \text{on}\ \ V_e(1/2).
$$
We call $\sigma_e$ the {\it second Miyamoto involution} associated 
to $e$.

\begin{rem}
  It is shown in \cite{C} and \cite{M1} that the Miyamoto involution $\tau_e$ 
  associated to a conformal vector $e$ of the moonshine VOA \cite{FLM} 
  with central charge 1/2 defines a 2A-involution of the Monster.
\end{rem}

\subsection{Commutant superalgebra}

We keep the same notation as in the previous subsection.
Let $V$ be a simple VOA of CFT-type and $e\in V$ a rational conformal
vector with central charge 1/2.
Set $T_e(h):= \{ v\in V \mid L^e(0) v= h\cd v\}$ for $h=0,1/2,1/16$.
Then $T_e(h)$ is a space of highest weight vectors for $\vir (e)$ and is 
canonically isomorphic to $\hom_{\vir (e)}(L(1/2,h), V)$ for $h=0,1/2,1/16$.
Therefore, we have a decomposition as follows:
$$
  V=L(1/2,0)\tensor T_e(0)\oplus L(1/2,1/2)\tensor T_e(1/2)\oplus 
    L(1/2,1/16)\tensor T_e(1/16).
$$

\begin{lem}
  A decomposition $\w=e+(\w-e)$ is orthogonal.
\end{lem}

\pf
We compute $e_{(1)} \w_{(2)} e$. 
$$
\begin{array}{ll}
  e_{(1)} \w_{(2)} e 
  &= \w_{(2)} e_{(1)} e+[e_{(1)},\w_{(2)}]e
  \vsb\\
  &= 2\w_{(2)} e -[\w_{(2)},e_{(1)}]e
  \vsb\\
  &=  2\w_{(2)} e-\{ (\w_{(0)} e)_{(3)}+2(\w_{(1)} e)_{(2)} 
      +(\w_{(2)}e)_{(1)}\} e
  \vsb\\
  &= 2\w_{(2)} e -(\w_{(2)}e)_{(1)}e .
\end{array}
$$
By the skew-symmetry, we have 
$
  (\w_{(2)} e)_{(1)} e= e_{(1)} \w_{(2)} e - \w_{(0)} e_{(2)} \w_{(2)} e.
$
Since $e_{(2)} \w_{(2)} e \in V_0=\C \vac$, $\w_{(0)} e_{(2)} \w_{(2)} e=0$ 
and so $(\w_{(2)} e)_{(1)} e= e_{(1)}\w_{(2)} e$.
Substituting this into the equality above, we get 
$e_{(1)}\w_{(2)} e = \w_{(2)} e$.
Namely, $\w_{(2)} e$ is an eigenvector for $e_{(1)}$ with eigenvalue 1.
Since $V$ is a module for $\vir (e)$, there is no eigenvector with 
$e_{(1)}$-weight 1. Hence $\w_{(2)} e=0$. 
Then the assertion follows from Theorem 5.1 of \cite{FZ}.
\qed
\vsb

Recall the commutant subalgebra $\com_V(\vir(e)):=\ker_V e_{(0)}$ 
defined in \cite{FZ}.
By the lemma above, $(T_e(0),\w-e)$ forms a sub VOA of $V$ whose action
on $V$ is commutative with that of $\vir (e)$ on $V$.
In particular, $T_e(h)$, $h=0,1/2,1/16$, are $T_e(0)$-modules.

\begin{prop}
  (1) $T_e(0)=\ker_V e_{(0)} =\com_V(\vir(e))$ is a simple sub 
  VOA with the Virasoro vector $\w-e$.\\
  (2) $T_e(1/2)$ is an irreducible $T_e(0)$-module.
  \\
  (3) $\vir (e)=\ker_V (\w-e)_{(0)}=\com_V (T_e(0))$.
\end{prop}

\pf
(1): 
Let $v\in V$.
Since $e_{(1)} v=0$ implies $e_{(0)} v=0$, $T_e(0)=\ker_V e_{(1)}
=\ker_V e_{(0)}$.
So we only need to show that $T_e(0)$ is simple.
Since $V$ is simple, the $\tau_e$-orbifold $V^{\la \tau_e\ra}
=V_e(0)\oplus V_e(1/2)$ is simple.
Then the $\sigma_e$-orbifold 
$(V^{\la \tau_e\ra})^{\la \sigma_e\ra}=V_e(0)$ is also simple.
Since $\vir (e)\tensor T_e(0)\ni a\tensor b\mapsto a_{(-1)}b\in V_e(0)$ 
is an isomorphism of VOAs, $T_e(0)$ is also simple.
\\
(2): Since both $V^{\la \tau_e \ra}= V_e(0)\oplus V_e(1/2)$ and
$V_e(0)$ are simple VOAs, $V_e(1/2)$ is an irreducible 
$V_e(0)$-module. So $T_e(1/2)$ is also irreducible.
\\ 
(3): As $\w-e$ is a conformal vector,
$\ker_V (\w-e)_{(0)}$ is generally contained in $\ker_V (\w-e)_{(1)}$.
On the other hand, since $V$ is of CFT-type, $\ker_V(\w-e)_{(1)}=\vir(e)$.
Then
$$
  \vir (e)\subset \com_V(\com_V(\vir(e))=\ker_V(\w-e)_{(0)}
$$
implies $\vir (e)=\ker_V (\w -e)_{(0)}$.
\qed

\begin{thm}\label{SVOA}
  Suppose that $T_e(1/2)\ne 0$.
  Then there exists a simple SVOA structure on $T_e(0)\oplus T_e(1/2)$
  such that the even part of a tensor product of SVOAs
  $$
    \{ L(1/2,0)\oplus L(1/2,1/2)\} \tensor 
    \{ T_e(0)\oplus T_e(1/2) \}
  $$ 
  is isomorphic to $V_e(0)\oplus V_e(1/2)$ as a VOA.
\end{thm}

\pf
We shall define vertex operators on an abstract space $T_e(0)
\oplus T_e(1/2)$.
First, we show an existence of a $T_e(0)$-intertwining operator
of type $T_e(1/2)\times T_e(1/2)\to T_e(0)$.
Write $Y_V(\w,z)=\sum_{n\in \Z} L(n)z^{-n-2}$ and $Y_V(e,z)=\sum_{n\in \Z}
L^e(n)z^{-n-2}$.
Since $L(0)-L^e(0)$ semisimply acts on both $T_e(0)$ and $T_e(1/2)$, 
we can take bases $\{ a^\gamma \mid \gamma\in \Gamma\}$ and 
$\{ u^\lambda \mid \lambda \in \Lambda\}$ of $T_e(0)$ and 
$T_e(1/2)$, respectively, consisting of eigenvectors for
$L(0)-L^e(0)$.
Let $\pi_\gamma : V_e(0)\to L(1/2,0)\tensor a^\gamma$, 
$\gamma \in \Gamma$, be a projection map.
For $\gamma \in \Gamma$ and $\lambda,\mu\in \Lambda$, 
we define a linear operator $I_{\lambda\mu}^\gamma(\cd,z)$
of type $L(1/2,1/2) \times L(1/2,1/2) \to L(1/2,0) 
\tensor a^\gamma$ by
$$
\begin{array}{ll}
  I_{\lambda\mu}^\gamma(x,z)y
  &:= z^{-L(0)+L^e(0)} \pi_\gamma  
  Y(z^{L(0)-L^e(0)} x\tensor u^\lambda,z) z^{L(0)-L^e(0)} y\tensor u^\mu
  \vsb\\
  &= z^{-\abs{\gamma}+\abs{\lambda}+\abs{\mu}}
    \pi_\gamma Y(x\tensor u^\lambda,z) y\tensor u^\mu,
\end{array}
$$
for $x,y\in L(1/2,1/2)$, where $\abs{\lambda}$, $\abs{\mu}$ and 
$\abs{\gamma}$ denote the $(L(0)-L^e(0))$-weight of $u^\lambda$, $u^\mu$ 
and $a^\gamma$, respectively.
Then by \cite{DL} \cite{M1} 
the operator $I_{\lambda\mu}^\gamma(\cd,z)$ is an $L(1/2,0)$-intertwining 
operator of type $L(1/2,1/2)\times L(1/2,1/2)\to L(1/2,0)$.
Since the space of intertwining operators of that type is one-dimensional, 
each $I_{\lambda \mu}^{\gamma}(\cd,z)$ is proportional to 
the vertex operator map $Y_M(\cd,z)$ on the SVOA 
$M=L(1/2,0) \oplus L(1/2,1/2)$ which we constructed explicitly in 
Section \ref{Ising}.
Thus there exist scalars $c_{\lambda\mu}^\gamma \in \C$ such that 
$I_{\lambda\mu}^\gamma (\cd,z) = c_{\lambda\mu}^\gamma Y_M(\cd,z)$.
Then the vertex operator of $x\tensor u^\lambda\in L(1/2,1/2)
\tensor T_e(1/2)$ on $V_e(1/2)$ can be written as follows:
$$
  Y_V(x\tensor u^\lambda,z) y\tensor u^\mu 
  = Y_M(x,z)y\tensor 
    \dsum_{\gamma \in \Gamma} c_{\lambda\mu}^\gamma a^\gamma 
    z^{\abs{\gamma}-\abs{\lambda}-\abs{\mu}}.
$$
Thus, by setting $J(u^\lambda,z)u^\mu:= \dsum_{\gamma\in \Gamma}
c^\gamma_{\lambda\mu} a^\gamma z^{\abs{\gamma}-\abs{\lambda}-
\abs{\mu}}$, we obtain a decomposition
$$
  Y_V(x\tensor u^\lambda,z) y\tensor u^\mu 
  =Y_M(x,z)y\tensor J(u^\lambda,z)u^\mu
$$
for $x\tensor u^\lambda$, $y\tensor u^\mu\in L(1/2,1/2)\tensor 
V_e(1/2)$.
We claim that $J(\cd,z)$ is a $T_e(0)$-intertwining operator of type 
$T_e(1/2)\times T_e(1/2)\to T_e(0)$.
It is obvious that $J(u,z)v$ contains finitely many negative 
powers of $z$ and the $(\w-e)_{(0)}$-derivation property 
$J((\w-e)_0u,z)v=\fr{d}{dz} J(u,z)v$ hold for all $u, v\in T_e(1/2)$.
So we should show that $J(\cd,z)$ satisfies both the commutativity and 
the associativity.
Let $a\in T_e(0)$ and $u,v\in T_e(1/2)$ be arbitrary elements.
Then the commutativity of vertex operators on $V$ gives 
$$
\begin{array}{l}
  (z_1-z_2)^N Y_V(\vac\tensor a,z_1) 
    Y_V(\psi_{-\hf}\vac \tensor u,z_2) 
    \psi_{-\hf}\vac \tensor v
  \vsb\\
  = (z_1-z_2)^N Y_V(\psi_{-\hf}\vac\tensor u,z_2)
    Y_V(\vac\tensor a,z_1) \psi_{-\hf}\vac\tensor v
\end{array}
$$
for sufficiently large $N$.
Rewriting the equality above we get 
$$
\begin{array}{l}
  (z_1-z_2)^N Y_M (\psi_{-\hf}\vac,z_2)
    \psi_{-\hf} \vac \tensor 
    Y_{T_e(0)}(a,z_1)J(u,z_2)v
  \vsb\\
  =(z_1-z_2)^N Y_M(\psi_{-\hf}\vac,z_2) \psi_{-\hf} \vac \tensor 
    J(u,z_2)Y_{T_e(1/2)}(a,z_1)v,
\end{array}
$$
where $Y_{T_e(0)}(a,z)$ and $Y_{T_e(1/2)}(\cd,z)$ denote the vertex operator 
of $a\in T_e(0)$ on $T_e(0)$ and $T_e(1/2)$, respectively.
By comparing the coefficients of $(\psi_{-\hf}\vac)_{(0)}
\psi_{-\hf}\vac=\vac$, we get the commutativity:
$$
  (z_1-z_2)^N Y_{T_e(0)} (a,z_1) J(u,z_2) v
  = (z_1-z_2)^N J(u,z_2) Y_{T_e(1/2)}(a,z_1)v.
$$
Similarly, by considering coefficients of
$Y_V(Y_V(\vac\tensor a,z_0)\psi_{-\hf}\vac\tensor u,z_2)
\psi_{-\hf}\vac\tensor v$ in $V$, we obtain the associativity:
$$
  (z_0+z_2)^N Y_{T_e(0)}(a,z_0+z_2) J(u,z_2) v
  = (z_2+z_0)^N J(Y_{T_e(1/2)}(a,z_0) u,z_2)v.
$$
Hence, $J(\cd,z)$ is a $T_e(0)$-intertwining operator of 
the desired type.

Using $Y_V(\cd,z)$ and $J(\cd,z)$, we introduce a vertex operator map
$\hat{Y}(\cd,z)$ on $T_e(0)\oplus T_e(1/2)$.
Let $a,b\in T_e(0)$ and $u,v\in T_e(1/2)$.
We define 
$$
\begin{array}{l}
  \vac\tensor \hat{Y}(a,z)b
    := Y_V(\vac\tensor a,z)\vac\tensor b,
  \q \psi_{-\hf}\vac \tensor \hat{Y}(a,z)u
    := Y_V(\vac\tensor a,z)\psi_{-\hf}\vac \tensor u,
  \vsb\\
  \psi_{-\hf}\vac \tensor \hat{Y}(u,z)a
    := e^{z(L(-1)-L^e(-1)} Y_V(\vac\tensor a,z) 
    \psi_{-\hf}\vac \tensor u,
  \q \hat{Y}(u,z)v:= J(u,z)v.
\end{array}
$$
Then all $\hat{Y}(\cd,z)$ are $T_e(0)$-intertwining operators.
We note that $\hat{Y}(\cd,z)$ satisfies the vacuum condition:
$$
  \hat{Y}(x,z)\vac\in x+\l( T_e(0)\oplus T_e(1/2)\r) [[z]]z
$$
for any $x\in T_e(0)\oplus T_e(1/2)$.
Hence, to prove that $T_e(0)\oplus T_e(1/2)$ is a simple SVOA, 
it is sufficient to show that the vertex operator map 
$\hat{Y}(\cd,z)$ defined above satisfies the commutativity.
By our definition, the vertex operator map $Y_V(a\tensor x,z)$ of 
$a\tensor x\in L(1/2,h)\tensor T_e(h)=V_e(h)$, $h=0,1/2$, can be 
written as $Y_M(a,z)\tensor \hat{Y}(x,z)$.
Because of our manifest construction of $Y_M(\cd,z)$ in Section \ref{Ising}, 
we can perform explicit computations of the vertex operator $Y_M(\cd,z)$ 
on $L(1/2,0)\oplus L(1/2,1/2)$.
Therefore, by comparing the coefficients of vertex operators on $V$, 
we can  prove that $\hat{Y}(\cd,z)$ satisfies the (super-)commutativity.
Thus, by our definition, $(T_e(0)\oplus T_e(1/2),\hat{Y}(\cd,z),
\vac,\w-e)$ carries a structure of a simple SVOA.
The rest of the assertion is now clear.
\qed

\begin{rem}
  There is another proof of Theorem \ref{SVOA} in \cite{Ho1}.
  In \cite{Ho1}, he assumed the existence of a positive definite 
  invariant bilinear form on a real form of $V$.
  However, our argument does not need the assumption on the 
  unitary form.
\end{rem}

Since $\tau_e^2=1$ on $V$, the space $V_e(1/16)$ is an irreducible 
$V^{\la \tau_e\ra}$-module.
As a $(V^{\la \tau_e\ra})^{\la \sigma_e\ra}=\vir (e) \tensor
T_e(0)$-module, $V_e(1/16)$ can be written as $L(1/2,1/16) \tensor 
T_e(1/16)$. 
It is not clear that $T_e(1/16)$ is irreducible under $T_e(0)$. 
However, we can prove that it is irreducible under $T_e(0)\oplus T_e(1/2)$.

\begin{thm}\label{twisted piece}
  Suppose that $V_e(1/2)\ne 0$ and $V_e(1/16)\ne 0$.
  Then $T_e(1/16)$ carries a structure of an irreducible
  $\Z_2$-twisted $T_e(0)\oplus T_e(1/2)$-module.
  Moreover, $V_e(1/16)$ is isomorphic to a tensor product of an 
  irreducible $\Z_2$-twisted $L(1/2,0)\oplus   L(1/2,1/2)$-module 
  $L(1/2,1/16)^+$ and an irreducible $\Z_2$-twisted $T_e(0)\oplus 
  T_e(1/2)$-module $T_e(1/16)$.
\end{thm}

\pf
The idea of the proof is the same as that of Theorem \ref{SVOA}.
Computing vertex operators on $L(1/2,1/16)^+$ and 
then comparing the coefficients in $V$, we will reach the assertion.
Denote by $Y_N(\cd,z)$ the vertex operator map on the $\Z_2$-twisted 
$L(1/2,0)\oplus L(1/2,1/2)$-module $L(1/2,1/16)^+$ as we constructed 
in Section \ref{Ising}.
Let $a\tensor b\in L(1/2,h)\tensor T_e(h)$ with $h=0$ or $1/2$
and $x\tensor y\in L(1/2,1/16)\tensor T_e(1/16)$.
As we did before, we can find $T_e(0)$-intertwining operators
$Y_{T_e(h)\times T_e(1/16)}(\cd,z)$ of types 
$T_e(h)\times T_e(1/16)\to T_e(1/16)$ such that 
\begin{equation}
  Y_V(a\tensor b,z)x\tensor y 
  = Y_N (a,z)x \tensor Y_{T_e(h)\times T_e(\fr{1}{16})}(b,z)y .
\end{equation}
Define $\hat{Y}(b,z)y:= Y_{T_e(h)\times T_e(\fr{1}{16})}(b,z)y$
for $b\in T_e(h)$, $h=0,1/2$ and $y\in T_e(1/16)$.
By direct computations, we can prove that the $\Z_2$-twisted 
Jacobi identity for $Y_N(\cd,z)$ together with the Jacobi identity for 
$Y_V(\cd,z)$ gives the $\Z_2$-twisted Jacobi identity for 
$\hat{Y}(\cd,z)$.
Thus, $(T_e(1/16),\hat{Y}(\cd,z))$ is a $\Z_2$-twisted 
$T_e(0)\oplus T_e(1/2)$-module.
Since $V_e(1/16)=L(1/2,1/16)\tensor T_e(1/16)$ is irreducible under 
$V_e(0)\oplus V_e(1/2)$, the irreducibility of $T_e(1/16)$ is obvious.
\qed

\subsection{One point stabilizer}

In the rest of this section we will work the following setup:

\begin{hypo}\label{II}
  \q \\
  (1) $V$ is a holomorphic VOA of CFT-type.
  \\
  (2) $e$ is a rational conformal vector of $V$ with central charge 1/2.
  \\
  (3) $V_e(h)\ne 0$ for $h=0,1/2,1/16$.
  \\
  (4) $V_e(0)$ and $T_e(0)$ are rational $C_2$-cofinite VOAs of CFT-type.
  \\
  (5) $V_e(1/16)$ is a simple current $V^{\la \tau_e\ra}=V_e(0)\oplus 
  V_e(1/2)$-module.
  \\
  (6) $T_e(1/2)$ is a simple current $T_e(0)$-module.
\end{hypo}

Define the one-point stabilizer by 
$C_{\aut (V)}(e):= \{ \rho \in \aut (V) \mid \rho (e)=e\}$.
Then by $\tau_{\rho(e)}=\rho \tau_e \rho^{-1}$ for any $\rho \in \aut (V)$,
we have $C_{\aut(V)}(e)\subset  C_{\aut(V)}(\tau_e)$, 
where $C_{\aut(V)}(\tau_e)$ denotes the centralizer of an involution 
$\tau_e\in \aut (V)$.

\begin{lem}\label{lem:8.1.7}
  There are group homomorphisms $\psi_1 :C_{\aut(V)}(e)\to 
  C_{\aut(V^{\la \tau_e\ra})}(e)$ and $\psi_2: C_{\aut(V^{\la \tau_e\ra})}(e) 
  \to \aut (T_e(0))$ such that $\ker (\psi_1)=\la \tau_e\ra$ and 
  $\ker(\psi_2)=\la \sigma_e\ra$.
\end{lem}

\pf
Let $\rho \in C_{\aut (V)}(e)$.
Then $\rho$ preserves the space of highest weight vectors $T_e(h)$ where 
$h\in \{ 0,1/2,1/16\}$.
Then we can define the actions of $\rho$ on the space of highest weight 
vectors $T_e(h)$ and the components $V_e(h)$ for $h \in \{ 0,1/2,1/16\}$.
In particular, we have group homomorphisms 
$\psi_1: C_{\aut (V)}(e)\to C_{\aut (V^{\la \tau_e\ra})}(e)$ and 
$\psi_2: C_{\aut (V^{\la \tau_e\ra})}(e)\to \aut (T_e(0))$ 
by a natural way.
Assume that $\psi_1(\rho)=\id_{V^{\la \tau_e\ra}}$ for $\rho \in 
C_{\aut (V)}(e)$.
Since $\rho\in C_{\aut (V)}(\tau_e)$, $\rho$ acts on $V_e(1/16)$ 
and commutes with the action of $V^{\la \tau_e\ra} =V_e(0)\oplus V_e(1/2)$ 
on its module $V_e(1/16)$.
Therefore, $\rho$ on $V_e(1/16)$ is a scalar by Schur's lemma and hence 
$\rho \in \la \tau_e\ra \subset C_{\aut (V)}(\tau_e)$.
Similarly, if $\psi_1(\rho')=\id_{T_e(0)}$ for $\rho'\in 
C_{\aut(V^{\la \tau_e\ra})}(e)$, then  $\rho'\in \la \sigma_e\ra
\subset C_{\aut(V^{\la \tau_e\ra})}(e)$.
\qed

\begin{thm}\label{thm:8.1.8}
  Under Hypothesis \ref{II}, 
  $V^{\la \tau_e\ra}$ has exactly four inequivalent 
  irreducible modules, $V^{\la \tau_e\ra}$, $V_e(1/16)$, 
  $W^0:=L(1/2,0)\tensor T_e(1/2)\oplus L(1/2,1/2)\tensor T_e(0)$ 
  and 
  $$
    W^1:=V_e(1/16)\fusion_{V^{\la \tau_e\ra}} W^0 .
  $$
\end{thm}

\pf
Note that $V_e(0)=\vir (e) \tensor T_e(0)$ and $V^{\la \tau_e\ra}$ are
simple rational $C_2$-cofinite VOAs of CFT-type under Hypothesis \ref{II}.
Therefore, we can apply a theory of fusion products here.
Since $V=V^{\la \tau_e\ra}\oplus V_e(1/16)$ is a $\Z_2$-graded simple 
current extension of $V^{\la \tau_e\ra}$, every irreducible 
$V^{\la \tau_e\ra}$-module is lifted to be either an irreducible $V$-module
or an irreducible $\tau_e$-twisted $V$-module.
Moreover, the $\tau_e$-twisted $V$-module is unique up to isomorphism 
by Theorem 10.3 of \cite{DLM2}.
Consider a $V_e(0)$-module $L(1/2,1/2)\tensor T_e(0)$.
Since $T_e(1/2)$ is a simple current $T_e(0)$-module, the space 
$$
  W^0= L(1/2,1/2)\tensor T_e(0)\oplus L(1/2,0)\tensor T_e(1/2)
$$
has a unique structure of an irreducible $V^{\la \tau_e\ra}$-module
by Theorem \ref{induced module}.
We note that the top weight of $W^0$ is half-integral.
Thus the induced module 
$$
  W=W^0\oplus W^1,\q W^1=V_e(1/16)\fusion_{V^{\la \tau_e\ra}} W^0,
$$
becomes an irreducible $\tau_e$-twisted $V$-module again by Theorem
\ref{induced module}.
Therefore, $V^{\la \tau_e\ra}$ has exactly four irreducible modules
as in the assertion.
Finally we remark that $V^{\la \tau_e\ra}$, $V_e(1/16)$ and $W^1$ have 
integral top weights.
\qed

By the fusion rules \eqref{fusion ising}, we note that $W^1$ as a 
$\vir(e)$-module is a direct sum of copies of $L(1/2,1/16)$.
Set the space of highest weight vectors of $W^1$ by 
$Q_e(1/16):=\{ v\in W^1 \mid L^e(0)v=(1/16)\cd v\}$.
Then as a $\vir(e) \tensor T_e(0)$-module, 
$W^1\simeq L(1/2,1/16)\tensor Q_e(1/16)$.
In this case, we can also verify that the space $Q_e(1/16)$ naturally 
carries an irreducible $\Z_2$-twisted $T_e(0)\oplus T_e(1/2)$-module 
structure.

\begin{prop}\label{prop:8.1.9}
  If the $\Z_2$-twisted $T_e(0)\oplus T_e(1/2)$-module $T_e(1/16)$ is 
  irreducible as a $T_e(0)$-module, then its $\Z_2$-conjugate is 
  isomorphic to $Q_e(1/16)$ as a $\Z_2$-twisted $T_e(0)\oplus T_e(1/2)$-module.
  In this case there are three irreducible $T_e(0)$-modules, $T_e(0)$, 
  $T_e(1/2)$ and $T_e(1/16)$.
  Conversely, if $T_e(1/16)$ as a $T_e(0)$-module is not irreducible, 
  then so is $Q_e(1/16)$ and in this case there are six inequivalent 
  irreducible $T_e(0)$-modules.
\end{prop}

\pf
Assume that $T_e(1/16)$ is irreducible as a $T_e(0)$-module. 
Then its $\Z_2$-conjugate is not isomorphic to $T_e(1/16)$ as a 
$\Z_2$-twisted $T_e(0)\oplus T_e(1/2)$-module.
We denote the $\Z_2$-conjugate of $T_e(1/16)$ by $T_e(1/16)^-$.
It is shown in Theorem \ref{super induced module} that every irreducible 
$T_e(0)$-module is lifted to be either 
an irreducible $T_e(0)\oplus T_e(1/2)$-module or an irreducible 
$\Z_2$-twisted $T_e(0)\oplus T_e(1/2)$-module.
Then by the classification of irreducible $V^{\la \tau_e\ra}$-modules in
Theorem \ref{thm:8.1.8}, we see that any $\Z_2$-twisted irreducible 
$T_e(0)\oplus T_e(1/2)$-module is isomorphic to one and only one of 
$T_e(1/16)$ and $T_e(1/16)^-=Q_e(1/16)$.

Conversely, if $T_e(1/16)$ is not irreducible, then it is a direct sum
of two inequivalent irreducible $T_e(0)$-module as $T_e(1/2)$ is a 
simple current $T_e(0)$-module.
Then $Q_e(1/16)$ is also a direct sum of two inequivalent irreducible 
$T_e(0)$-modules and $Q_e(1/16)\not\simeq T_e(1/16)$ as $T_e(0)$-modules
because of the classification of irreducible $V^{\la \tau_e\ra}$-modules.
\qed

\begin{cor}\label{cor:8.1.10}
  If $T_e(1/16)$ is irreducible as a $T_e(0)$-module, then 
  $V^{\la \tau_e\ra}\oplus W^1$ is a $\Z_2$-graded simple current
  extension of $V^{\la\tau_e\ra}$ which is isomorphic to 
  $V=V^{\la \tau_e\ra} \oplus V_e(1/16)$.
\end{cor}

\pf
If $T_e(1/16)$ is an irreducible $T_e(0)$-module, then by the 
previous proposition the $\sigma_e$-conjugate 
$V_e(0)\oplus V_e(1/2)$-module of 
$V_e(1/16)=L(1/2,1/16)\tensor T_e(1/16)$ is isomorphic to 
$W^1=L(1/2,1/16)\tensor Q_e(1/16)$.
Then as the $\sigma_e$-conjugate extension of $V^{\la\tau_e\ra}
\oplus V_e(1/16)$, $V^{\la \tau_e\ra}\oplus W^1$ has a structure of a 
$\Z_2$-graded extension.
\qed

\begin{rem}
  The above corollary implies that the $\Z_2$-twisted orbifold 
  construction applied to $V$ in the case of $\Z_2=\la \tau_e\ra$ yields
  again $V$ itself.
\end{rem}

\begin{thm}\label{one-point stabilizer}
  Under Hypothesis \ref{II}, 
  \\
  (i)\ $\psi_2$ is surjective, that is, 
  $C_{\aut (V^{\la \tau_e\ra})} (e)\simeq \la \sigma_e\ra .\aut (T_e(0))$.
  \\
  (ii) $\aut (T_e(0)\oplus T_e(1/2))\simeq 2. (C_{\aut(V^{\la \tau_e\ra})}
  (e)/\la \sigma_e\ra)$, where $2$ denotes the canonical $\Z_2$-symmetry on 
  the SVOA $T_e(0)\oplus T_e(1/2)$.
  \\
  (iii) $\abs{C_{(\aut (V^{\la \tau_e\ra}))}(e):C_{\aut(V)}(e)/\la \tau_e\ra} 
  \leq 2$.
  \\
  (iv) If $C_{\aut (V)}(e)/\la \tau_e\ra$ is simple or has an odd order, 
  then extensions in (i) and (ii) split.
  That is, $C_{\aut(V^{\la \tau_e\ra})}(e) \simeq \la \sigma_e\ra 
  \times C_{\aut (V)}(e)/\la \tau_e\ra$ and $\aut (T_e(0)\oplus T_e(1/2))
  \simeq 2\times \aut (T_e(0))$.
\end{thm}

\pf
We have an injection from $C_{\aut (V^{\la \tau_e\ra})}(e)/\la 
\sigma_e\ra$ to $\aut (T_e(0))$ by Lemma \ref{lem:8.1.7}. 
We show that every element in $\aut (T_e(0))$ lifts to be an element
in $C_{\aut(V^{\la \tau_e\ra})}(e)$.
By Proposition \ref{prop:8.1.9}, every irreducible $T_e(0)$-module appears 
in one of $T_e(0)$, $T_e(1/2)$, $T_e(1/16)$ or $Q_e(1/16)$ as a submodule.
In particular, we find that $T_e(0)$ is the only irreducible $T_e(0)$-module 
whose top weight is integral and $T_e(1/2)$ is the only irreducible 
$T_e(0)$-module whose top weight is in $1/2+\N$.
Let $\rho \in \aut (T_e(0))$.
Then by considering top weights we can immediately see that 
$T_e(0)^\rho \simeq T_e(0)$ and $T_e(1/2)^\rho\simeq T_e(1/2)$.
Then by Theorem \ref{Sh} we have a lifting $\tilde{\rho}\in 
\aut (T_e(0)\oplus T_e(1/2))$ such that $\tilde{\rho}T_e(0)=T_e(0)$,
$\tilde{\rho}T_e(1/2)=T_e(1/2)$ and $\tilde{\rho}|_{T_e(0)}=\rho$.
Since this lifting is unique up to a multiple of the canonical 
$\Z_2$-symmetry on $T_e(0)\oplus T_e(1/2)$, we have 
$\aut (T_e(0)\oplus T_e(1/2))\simeq 2.\aut (T_e(0))$.
Now consider the canonical extension of $\tilde{\rho}$ to 
$C_{\aut (V^{\la \tau_e\ra})}(e)$. 
We define $\tilde{\tilde{\rho}}\in C_{\aut(V^{\la \tau_e\ra})}(e)$ by
$$
  \tilde{\tilde{\rho}}|_{L(1/2,h)\tensor T_e(h)}
  =\id_{L(1/2,h)}\tensor \tilde{\rho}
$$
for $h=0,1/2$.
Then by this lifting $C_{\aut (V^{\la \tau_e\ra})}(e)$ contains 
a subgroup to isomorphic to $2.\aut (T_e(0))$.
Moreover, the canonical $\Z_2$-symmetry on the SVOA $T_e(0)\oplus T_e(1/2)$
is naturally extended to $\sigma_e\in C_{\aut (V^{\la \tau_e\ra})}(e)$.
Clearly $\psi_2(\tilde{\tilde{\rho}})=\rho$ and so $\psi_2$ is surjective.
Hence we have the desired isomorphisms $C_{\aut(V^{\la \tau_e\ra})}(e)
\simeq \la \sigma_e\ra. \aut (T_e(0))$ and 
$\aut (T_e(0)\oplus T_e(1/2))\simeq 2.
(C_{\aut(V^{\la \tau_e\ra})}(e)/\la \sigma_e\ra )$.
This completes the proof of (i) and (ii).

Consider (iii).
By Theorem \ref{thm:8.1.8}, there are exactly three irreducible 
$V^{\la \tau_e\ra}$-modules whose top weights are integral, namely,
$V^{\la \tau_e\ra}$, $V_e(1/16)$ and $W^1$.
Since $C_{\aut(V^{\la \tau_e\ra})}(e)$ acts on the 2-point set 
$\{ V_e(1/16), W^1\}$ as a permutation, there is a subgroup $H$ of 
$C_{\aut(V^{\la \tau_e\ra})}(e)$ of index at most 2 such that 
$V_e(1/16)^\pi \simeq V_e(1/16)$ as a $V^{\la \tau_e\ra}$-module for 
all $\pi \in H$.
Then there is a lifting $\tilde{\pi}\in C_{\aut (V)}(e)$ of $\pi$ such 
that $\psi_1(\tilde{\pi})=\pi$ for each $\pi \in H$ by Theorem \ref{Sh}.
Thus $\abs{C_{\aut(V^{\la \tau_e\ra})}(e):C_{\aut(V)}(e)/\la \tau_e\ra}
\leq 2$ and (iii) holds.

Consider (iv).
Suppose that $C_{\aut (V)}(e)/\la \tau_e\ra$ is a simple group or an odd 
group.
Then $C_{\aut(V^{\la \tau_e\ra})}(e)$ contains a simple 
group $C_{\aut (V)}(e)/\la \tau_e\ra$ with index at most 2 by (iii).
However, since $C_{\aut (V^{\la \tau_e\ra})}(e)$ contains a normal subgroup 
$\la \sigma_e\ra$ of order 2, the index $\abs{C_{\aut(V^{\la \tau_e\ra})}(e):
C_{\aut(V)}(e)/\la \tau_e\ra}$ must be 2 and hence we obtain the desired 
isomorphism $C_{\aut(V^{\la \tau_e\ra})}(e)=\la \sigma_e\ra \times 
C_{\aut(V)}(e)/\la \tau_e\ra$.
In this case, it is easy to see that the extension $\aut (T_e(0)\oplus 
T_e(1/2))=2.\aut(T_e(0))$ splits.
\qed

\begin{cor}\label{cor:3.15}
  If $C_{\aut (V)}(e)/ \la \tau_e\ra$ is simple, then $V_e(1/16)$ 
  is an irreducible $V_e(0)$-module and $T_e(1/16)$ is an irreducible 
  $T_e(0)$-module.
  Therefore, $V^{\la \tau_e\ra} \oplus W^1$ forms the $\sigma_e$-conjugate 
  extension of $V=V^{\la \tau_e\ra}\oplus V_e(1/16)$ and is isomorphic 
  to $V$.
\end{cor}

\pf
Let $H$ be the subgroup of $C_{\aut(V^{\la \tau_e\ra})}(e)$ which
fixes $V_e(1/16)$ in the action on the 2-point set $\{ V_e(1/16), W^1\}$.
It is shown in the proof of (iii) of Theorem \ref{one-point stabilizer} that 
we have inclusions 
$$
  H \subset C_{\aut (V)}(e)/\la \tau_e\ra 
  \subset C_{\aut (V^{\la \tau_e\ra})}(e) 
  = \la \sigma_e \ra \times C_{\aut (V)}(e)/\la \tau_e\ra .
$$
Therefore, $\sigma_e \not \in H$ and hence the $\sigma_e$ permutes $V_e(1/16)$ 
and $W^1$.
Then $V_e(1/16)$ is an irreducible $V_e(0)$-module by Proposition 
\ref{prop:8.1.9} and hence $T_e(1/16)$ as a $T_e(0)$-module is irreducible.
The rest of the assertion is now clear.
\qed

\begin{rem}
  A result similar to the assertion (iii) of Theorem 
  \ref{one-point stabilizer} is already established in \cite{M6}.
  Also, we should note that the idea of the above proof is already 
  developed in \cite{Sh}.
\end{rem}

\section{2A-framed VOA}

\begin{df} 
  A simple vertex operator algebra  $(V,\w)$ is called \emph{2A-framed}
  if there is an orthogonal decomposition $\w=e^1+\cds +e^n$ such that 
  each $e^i$ generates a sub VOA isomorphic to $L(1/2,0)$.
  The decomposition $\w=e^1+\cds +e^n$ is called a \emph{2A-frame} of 
  $V$.
\end{df}

\begin{rem}
  As shown in \cite{DMZ}, the Leech lattice VOA $V_{\Lambda}$ and
  the moonshine VOA $V^\natural$ are examples of 2A-framed VOAs.
\end{rem}

Let $(V,\w)$ be a 2A-framed VOA with a 2A-frame $\w=e^1+\cds +e^n$.
Set $T:=\vir(e^1)\tensor \cds \tensor \vir (e^n)$, where $\vir (e^i)$ 
denotes the sub VOA generated by $e^i$.
Then $T\simeq L(1/2,0)^{\tensor n}$ and $V$ is a direct sum of irreducible 
$T$-submodules $\tensor_{i=1}^n L(1/2,h_i)$ with $h_i\in \{ 0,1/2,1/16\}$. 
For each irreducible $T$-module $\tensor_{i=1}^n L(1/2,h_i)$, we associate 
its $1/16$-word $(\alpha_1,\cds,\alpha_n)\in \Z_2^n$ by the rule 
$\alpha_i=1$ if and only if $h_i=1/16$.
For each $\alpha \in \Z_2^n$, denote by $V^\alpha$ the sum of all 
irreducible $T$-submodules whose $1/16$-words are equal to $\alpha$ and 
define a linear code $S\subset \Z_2^n$ by $S=\{ \alpha \in \Z_2^n 
\mid V^\alpha\ne 0\}$.
Then we have the {\it 1/16-word decomposition} 
$V=\oplus_{\alpha \in S} V^\alpha$ of $V$.
By the fusion rules \eqref{fusion ising}, we have an $S$-graded structure 
$V^\alpha\cd V^\beta
\subset V^{\alpha +\beta}$.
Namely, the dual group $S^*$ of an abelian 2-group $S$ acts 
on $V$, and we find that this automorphism group coincides with 
the elementary abelian 2-group generated by the first Miyamoto involutions 
$\{ \tau_{e^i} \mid 1\leq i \leq n\}$.
Therefore, all $V^\alpha$, $\alpha \in S$, are irreducible 
$V^{S^*}=V^0$-modules by \cite{DM1}.
Since there is no $L(1/2,1/16)$-component in $V^0$, the fixed point subalgebra 
$V^{S^*}=V^0$ has the following shape:
$$
  V^0=\bigoplus_{h_i\in \{ 0,1/2\}} m_{h_1,\dots,h_n}
  L(1/2,h_1)\tensor \cds \tensor L(1/2,h_n),
$$
where $m_{h_1,\dots,h_n}$ denotes the multiplicity.
On $V^0$ we can define the second Miyamoto involutions $\sigma_{e^i}$ for 
$i=1,\dots,n$.
Denote by $Q$ the elementary abelian 2-subgroup of $\aut (V^0)$ generated 
by $\{ \sigma_{e^i} \mid 1\leq i\leq n\}$.
Then by the quantum Galois theory \cite{DM1} we have $(V^0)^Q=T$ and 
each $m_{h_1,\cds,h_n} L(1/2,h_1)\tensor \cds \tensor L(1/2,h_n)$ is 
an irreducible $T$-submodule.
Thus $m_{h_1,\cds,h_n}\in \{ 0,1\}$ and we obtain an even linear code 
$D:= \{ (2h_1,\cds ,2h_n)\in \Z_2^n \mid m_{h_1,\cds,h_n}\ne 0\}$ 
such that 
\begin{equation}\label{7.3.1}
  V^0=\bigoplus_{\alpha =(\alpha_1,\cds,\alpha_n)\in D} L(1/2,\alpha_1/2)
  \tensor \cds \tensor L(1/2,\alpha_n/2).
\end{equation}
We call a pair $(D,S)$ the {\it structure codes} of a 2A-framed VOA $V$.
Since powers of $z$ in an $L(1/2,0)$-intertwining operator of type 
$L(1/2,1/2)\times L(1/2,1/2)\to L(1/2,1/16)$ are half-integral, structure 
codes satisfy $D\subset S^\perp$.
Since all $V^\alpha$, $\alpha \in S$, are irreducible modules for $V^0$, 
the representation theory of $V^0$ is important to study a 2A-framed VOA $V$.
We review Miyamoto's results on the code VOAs \cite{M3} \cite{M4} in the 
next subsection.

\subsection{Code VOA}

Below we often identify the code $\Z_2^n$ with the power set of an 
$n$-point set $\Omega=\{ 1,2,\dots,n\}$ with the symmetric difference 
operation.
For $\alpha =(\alpha_1,\dots,\alpha_n)\in \Z_2^n$, we denote by 
$\supp (\alpha)$ the subset $\{ i \mid \alpha_i\ne 0\}$ of $\Omega$.
Let $D$ be a subcode of $\Z_2^n$.
Set $D^{(0)}:=\{ \alpha \in D \mid \la \alpha,\alpha \ra=0\}$ 
and $D^{(1)}:=\{ \alpha \in D \mid \la \alpha,\alpha \ra=1\}$.
For a code-word $\alpha =(\alpha_1,\dots,\alpha_n) \in D$, we define
$$
  U^\alpha:=L(1/2,\alpha_1/2)\tensor \cds \tensor L(1/2,\alpha_n/2).
$$

\begin{thm}\label{code VOA}
  (\cite{M2})
  For any linear code $D\subset \Z_2^n$, there exists a unique simple vertex 
  operator superalgebra structure on $U_D:=\oplus_{\alpha\in D} U^\alpha$ 
  as an extension of $U^0=L(1/2,0)^{\tensor n}$.
  The even part is $U_{D^{(0)}}:=\oplus_{\alpha \in D^{(0)}} U^\alpha$
  and is a $D^{(0)}$-graded simple current extension of $U^0$, and 
  the odd part is $U_{D^{(1)}}:=\oplus_{\beta \in D^{(1)}} U^\beta$.
\end{thm}

Let $D$ be an even subcode of $\Z_2^n$.
The simple VOA $U_D$ defined in the theorem above is called the 
{\it code VOA associated to a code $D$}.
The representation theory of $U_D$ is deeply studied in \cite{M3}.
We recall some results from \cite{M3}.
Since $U_D$ is a $D$-graded simple current extension of 
a rational VOA $U^0=L(1/2,0)^{\tensor n}$, it is also rational.
Let $M$ be an irreducible $U_D$-module.
Take an irreducible $U^0$-submodule $W$ of $M$.
Then $W$ is isomorphic to $\tensor_{i=1}^n L(1/2,h_i)$ with $h_i\in 
\{ 0,1/2,1/16\}$.
Define the 1/16-word $\tau(W)=(\alpha_1,\dots,\alpha_n) \in \Z_2^n$ of $W$ 
by $\alpha_i=1$ if and only if $h_i=1/16$.
Then by the fusion rules \eqref{fusion ising}, $\tau(W)$ is determined 
independently of a choice of $W$ and hence we can define the 1/16-word
$\tau(M)\in \Z_2^n$ of $M$ by $\tau(M):=\tau(W)$.

Assume that $\tau(M)=(0^n)$.
Then $W$ is isomorphic to $\tensor_{i=1}^n L(1/2,h_i)$ with $h_i\in \{ 
0,1/2\}$.
We set $\gamma:=(2h_1,\dots,2h_n)\in \Z_2^n$.
Since $L(1/2,0)$ is a simple current $L(1/2,0)$-module, we have 
$D_W=0$ and hence $M$ is uniquely determined by $W$ by Theorem
\ref{induced module} and has a shape
$$
  U_{D+\gamma}:= \oplus_{\alpha \in D+\gamma} U^{\alpha}.
$$
We call an irreducible $U_D$-module $U_{D+\gamma}$ a {\it coset module}.
By Theorem \ref{lifting} and the fusion rules \eqref{fusion ising}, 
we have the following fusion rules for $\gamma_1,\gamma_2\in \Z_2^n$:
\begin{equation}\label{fusion coset}
  U_{D+\gamma_1}\times U_{D+\gamma_2}=U_{D+\gamma_3}.
\end{equation}
In particular, $U_{D+\delta}\times U_{D+\delta}=U_D$ for any $\delta
\in \Z_2^n$ and hence all the coset modules are simple current
$U_D$-modules by Lemma \ref{criterion}.

\subsection{The Hamming code VOA}

Let $H_8$ be the $[8,4,4]$-Hamming code:
$$
  H_8:=\Span_{\Z_2}\{ (11111111), (11110000), (11001100), (10101010)\} .
$$
It is well-known that $H_8$ is the unique doubly even self-dual linear code 
of length 8 up to isomorphism.
Let us consider the Hamming code VOA $U_{H_8}$.
In order to construct 2A-framed VOAs, we will need some special properties 
that the Hamming code VOA $U_{H_8}$ possesses.
Roughly speaking, we can identify $L(1/2,1/16)$ with $L(1/2,0)$ and 
$L(1/2,1/2)$ by the symmetry of the Hamming code VOA.

Let $X$ be an irreducible $U_{H_8}$-module whose top weight is in $\hf\N$.
Then $\tau(X)=(0^8)$ or $(1^8)$ and if $\tau(X)=(0^8)$ then $X\simeq 
U_{H_8+\gamma}$ for some $\gamma \in \Z_2^8$.
If $\tau(X)=(1^8)$, then it is shown in \cite{M3} that there is a unique 
linear character $\chi$ on $H_8$ such that 
\begin{equation}
  X \simeq \ind_{U^0}^{U_{H_8}}(L(1/2,1/16)^{\tensor 8},\chi)
  = L(1/2,1/16)^{\tensor 8}\tensor_\C v_\chi,
\end{equation}
where $U^\alpha$, $\alpha \in H_8$, acts on $L(1/2,16)^{\tensor 8}$ 
by the fusion rule \eqref{fusion ising} and on $\C v_\chi$ as a scalar 
$\chi(\alpha)$ according to the character $\chi$.
Since the dual group $H_8^*$ of $H_8$ is naturally isomorphic to 
$\Z_2^8/H_8$, we can find a unique coset $\delta_\chi+H_8 \in \Z_2^8/H_8$ 
such that $\chi(\alpha)=\la \delta_\chi,\alpha\ra$ for all $\alpha \in H_8$.
So in the following we regard $\chi$ as an element in $\Z_2^8$.
Set $H(1/16,\chi)=L(1/2,1/16)^{\tensor 8} \tensor_\C v_\chi$ for 
$\chi \in \Z_2^8$.
Then $H(1/16,\chi_1)\simeq H(1/16,\chi_2)$ as $U_{H_8}$-modules if and 
only if $\chi_1-\chi_2\in H_8$ and the set of inequivalent irreducible 
$U_{H_8}$-modules whose top weights are contained in $\hf \N$ is given by 
$$
  \{ U_{H_8+\gamma},\ H(1/16,\chi) \mid \gamma +H_8, \chi+H_8 
  \in \Z_2^8/H_8\} .
$$
Surprisingly, we can identify a non-simple current 
$L(1/2,0)^{\tensor 8}$-module $L(1/2,1/16)^{\tensor 8}$ with a coset module 
as follows:

\begin{thm}\label{M4a}
  (\cite{M4})
  For each $\chi\in \Z_2^8$, there is an automorphism 
  $\sigma \in \aut (U_{H_8})$ such that the $\sigma$-conjugate 
  module $H(1/16,\chi)^\sigma$ is isomorphic to a coset module 
  $U_{H_8+\gamma}$ for some $\gamma\in \Z_2^8$ with 
  $\la \gamma,\gamma \ra =1$.
  In particular, $H(1/16,\chi)$ is a simple current $U_{H_8}$-module.
\end{thm}

\begin{cor}\label{M4b}
  (\cite{M4})
  As a $\Z_2$-graded simple current extension of $U_{H_8}$, 
  there is a unique simple SVOA structure on $U_{H_8}\oplus 
  H(1/16,\chi)$ for all $\chi \in \Z_2^8$. 
\end{cor}

\pf
We can take an irreducible coset $U_{H_8}$-module $U_{H_8+\gamma}$ with 
$\la \gamma,\gamma\ra =1$ such that there is an automorphism 
$\sigma \in \aut (U_{H_8})$ such that the conjugate module 
$(U_{H_8+\gamma})^\sigma$ is isomorphic to $H(1/16,\chi)$ by 
Theorem \ref{M4a}.
Then $U_{H_8}\oplus U_{H_8+\gamma}$ and $U_{H_8}\oplus H(1/16,\chi)$ 
form mutually conjugate $\Z_2$-graded simple current extensions of 
$U_{H_8}$ under $\sigma \in \aut (U_{H_8})$.
Since $H_8\cup (H_8+\gamma)$ is an odd code, 
$U_{H_8}\oplus U_{H_8+\gamma}$ is a simple SVOA.
Then so is $U_{H_8}\oplus H(1/16,\chi)$.
\qed
\vsb

As an application of Proposition \ref{M4a}, the following fusion rules
are established in \cite{M4}:

\begin{thm}\label{M4c}
  (\cite{M4})
  We have the following fusion rules:
  $$
  \begin{array}{l}
    U_{H_8+\alpha} \times U_{H_8+\beta}= U_{H_8+\alpha+\beta}, 
    \vsb\\
    U_{H_8+\alpha} \times H(1/16,\beta)= H(1/16,\alpha+\beta), 
    \vsb\\
    H(1/16,\alpha) \times H(1/16,\beta)= U_{H_8+\alpha+\beta}, 
  \end{array}
  $$
  where $\alpha,\beta\in \Z_2^8$.
\end{thm}

Thanks to Corollary \ref{M4b} and Theorem \ref{M4c}, if an even linear 
code $D$ contains many subcodes isomorphic to the Hamming code $H_8$, 
then we can construct simple current extensions of the code VOA $U_D$
by using Theorem \ref{extension}.

\subsection{Construction of 2A-framed VOA}

In this subsection we give a refinement Miyamoto's construction of 2A-framed 
VOAs in \cite{M4}.
Here we assume the following:
\vsb

\begin{hypo}\label{I}
\q
\\
(1) $(D,S)$ is a pair of even linear even codes of $\Z_2^n$ 
such that 
\vsb\\
\begin{tabular}{p{440pt}}
  (1-i)\ $D\subset S^\perp$,
  \vsb\\ 
  (1-ii)\ for each $\alpha \in S$, there is a subcode $E^\alpha\subset D$ 
  such that $E^\alpha$ is a direct sum of the Hamming code $H_8$ and 
  $\supp (E^\alpha)=\supp (\alpha)$, where $\supp (A)$ denotes 
  $\cup_{\beta \in A}\supp(\beta)$ for a subset $A$ of $\Z_2^n$.
\end{tabular}
\vsb\\
(2) $V^0=U_D$ is the code VOA associated to the code $D$.
\vsb\\
(3) $\{ V^\alpha \mid \alpha \in S\}$ is a set of irreducible $V^0$-modules 
such that 
\vsb\\
\begin{tabular}{p{440pt}}
  (3-i)\ $\tau(V^\alpha)=\alpha$ for all $\alpha \in S$, 
  \vsb\\
  (3-ii)\ all $V^\alpha$, $\alpha \in S$, have integral top weights,
  \vsb\\
  (3-iii)\ the fusion product $V^\alpha \fusion_{V^0} V^\beta$ contains 
  at least one $V^{\alpha+\beta}$.
  That is, there is a non-trivial $V^0$-intertwining operator of type 
  $V^\alpha \times V^\beta \to V^{\alpha +\beta}$ for any $\alpha,\beta 
  \in S$.
\end{tabular}
\end{hypo}

Under  Hypothesis \ref{I} we will prove that $V:=\oplus_{\alpha \in S} V^\alpha$
has a structure of an $S$-graded simple current extension of $V^0$.
Before we begin the proof, we prepare some lemmas.

\begin{lem}\label{lem:7.4.3}
  Under Hypothesis \ref{I}, all $V^\alpha$, $\alpha \in S$, are 
  simple current $V^0$-modules and we have the fusion rules 
  $V^\alpha\times V^\beta=V^{\alpha +\beta}$ of $V^0$-modules 
   for all $\alpha,\beta\in S$.
\end{lem}

\pf
Suppose the fusion rule $V^\alpha \times V^\alpha=V^0$ of $V^0$-modules 
holds.
Then by Lemma \ref{criterion}, $V^\alpha$ is a simple current $V^0$-module 
because $V^0=U_D$ is a rational $C_2$-cofinite VOA of CFT-type.
Then by Hypothesis \ref{I} (3-iii) we have the desired fusion rule 
$V^\alpha \times V^\beta=V^{\alpha+\beta}$.
Therefore, we only prove the fusion rule $V^\alpha \times V^\alpha =V^0$ 
for each $\alpha \in S$.
By Hypothesis \ref{I} (1-i), $D$ contains a subcode $E^\alpha$ which is 
isomorphic to a direct sum of $H_8$ and $\supp (E^\alpha)=\supp (\alpha)$.
We may assume that $\alpha=(1^{8s}0^t)$ with $8s+t=n$.
Then $U_D$ contains a sub VOA 
$$
  L:=U_{E^\alpha}\tensor L(1/2,0)^{\tensor t}
  \simeq (U_{H_8})^{\tensor s}\tensor L(1/2,0)^{\tensor t}
$$ 
and $V^\alpha$ contains an irreducible $L$-submodule $X$ isomorphic to
$$
  H(1/16,\chi_1)\tensor \cds H(1/16,\chi_s)\tensor L(1/2,h_1)\tensor 
  \cds \tensor L(1/2,h_t)
$$
with $\chi_i\in \Z_2^8$, $1\leq i\leq s$, and $h_j\in \{ 0,1/2\}$, 
$1\leq j\leq t$.
Let $D=\sqcup_{i=0}^k (E^\alpha+\beta_i)$ be a coset decomposition.
We write $\beta_i=\gamma_i+\delta_i$ such that $\supp(\gamma_i)\subset 
\supp (\alpha)$ and $\supp (\delta_i)\cap \supp(\alpha)=\emptyset$.
Then $U_{E^\alpha +\beta_i}$ is isomorphic to 
$$
  U_{E^\alpha+\gamma_i}\tensor L(1/2,(\delta_i)_{8s+1}/2)\tensor \cds 
  \tensor L(1/2,(\delta_i)_n/2)
$$
as an $L$-module and $U_D=\oplus_{i=1}^k U_{E^\alpha+\beta_i}$ is a 
$D/E^\alpha$-graded simple current extension of $L$.
Then by the fusion rules \eqref{fusion ising} and Theorem \ref{M4c}, 
$(U_{E^\alpha+\beta_i})\fusion_L X$ is an irreducible $L$-module and 
$(U_{E^\alpha+\beta_i})\fusion_L X \not\simeq (U_{E^\alpha+\beta_j})
\fusion_L X$ unless $i=j$.
Therefore, $V^\alpha$ as an $L$-module is isomorphic to 
$V^\alpha= \oplus_{i=1}^k (V_{E^\alpha+\beta_i})\fusion_L X$.
Namely, $V^\alpha$ is a $D/E^\alpha$-stable $U_D$-module.
Then by Theorem \ref{lifting} together with fusion rules \eqref{fusion ising}
and Theorem \ref{M4c} we have a fusion rule $V^\alpha\times 
V^\alpha = V^0$ of $U_D$-modules which is a lifting of the fusion rule 
$X\times X= L$ of $L$-modules.
\qed

\begin{lem}\label{lem:7.4.4}
  Under Hypothesis \ref{I}, the space $V^0\oplus V^\alpha$ forms a simple VOA 
  as a $\Z_2$-graded simple current extension of $V^0$ for each 
  $\alpha \in S\setminus 0$.
\end{lem}

\pf
Here we use the same notation as in the proof of Lemma \ref{lem:7.4.3}.
By the coset decomposition $D=\sqcup_{i=1}^k (E^\alpha+\beta_i)$, 
$V^0=U_D=\oplus_{i=1}^k U_{E^\alpha+\beta_i}$ is a 
$D/E^\alpha$-graded simple current extension of $L=U_{E^\alpha}\tensor 
L(1/2,0)^{\tensor t}\simeq (U_{H_8})^{\tensor s} 
\tensor L(1/2,0)^{\tensor t}$.
By the fusion rule $X\times X=L$ 
of $L$-modules, $X$ is a simple current $L$-module by Lemma \ref{criterion}.
Then by the associativity of fusion products (cf.\ \cite{H4}), 
$(U_{E^\alpha+\beta_i})\fusion_L X$ is also a simple current $L$-module.
Thus we obtain the set of inequivalent simple current $L$-modules 
$$
  \mathcal{S}=\{\, U_{E^\alpha+\beta_i},\ (U_{E^\alpha+\beta_j})\fusion_L X 
  \mid 1\leq i,j\leq k\, \} 
$$
with the following $\l((D/E^\alpha)\oplus \Z_2\r)$-graded fusion rules:
$$
\begin{array}{l}
  U_{E^\alpha+\beta_i}\times U_{E^\alpha+\beta_j}
  = U_{E^\alpha+\beta_i+\beta_j},
  \vsb\\
  U_{E^\alpha+\beta_i}\times (U_{E^\alpha+\beta_j}\fusion_L X)
  = (U_{E^\alpha+\beta_i+\beta_j})\fusion_L X,
  \vsb\\
  (U_{E^\alpha+\beta_i}\fusion_L X)\times (U_{E^\alpha+\beta_j}\fusion_L X)
  = U_{E^\alpha+\beta_i+\beta_j}.
\end{array}
$$
Since $U_D=\oplus_{i=1}^k U_{E^\alpha+\beta_i}$ has a structure of a 
$D/E^\alpha$-graded simple current extension of $L$ and 
$L \oplus X$ has a structure of a $\Z_2$-graded simple current extension of 
$L$ by Corollary \ref{M4b}, we can apply Theorem \ref{extension} to 
$\mathcal{S}$ and hence we obtain a $\l((D/E^\alpha)\oplus \Z_2\r)$-graded 
simple current extension 
$$
  \l\{ \oplus_{i=1}^k U_{E^\alpha+\beta_i}\r\} \bigoplus 
  \l\{ \oplus_{i=1}^k (U_{E^\alpha+\beta_i})\fusion_L X \r\} 
$$
of $L$.
Since $V^0=\oplus_{i=1}^k U_{E^\alpha+\beta_i}$ and $V^\alpha = 
\oplus_{i=1}^k (U_{E^\alpha+\beta_i})\fusion_L X$, the $\Z_2$-graded space 
$V^0\oplus V^\alpha$ carries a simple VOA structure which is the desired
$\Z_2$-graded simple current extension of $V^0$.
\qed
\vsb

Now we can prove 

\begin{thm}\label{M4main}
  (\cite{M4})
  Under Hypothesis \ref{I}, the space $V=\oplus_{\alpha \in S} V^\alpha$ 
  has a unique structure of a simple VOA as an $S$-graded simple current 
  extension of $V^0$.
  In particular, there exists a 2A-framed VOA whose structure codes 
  are $(D,S)$.
\end{thm}

\pf
Let $\{ \alpha_1,\dots,\alpha_r\}$ be a linear basis of $S$ and 
set $S^i:=\Span_{\Z_2}\{ \alpha_1,\dots,\alpha_i\}$ for $1\leq i\leq r$.
We proceed by induction on $r$.
The case $r=0$ is trivial and the case $r=1$ is given by Lemma 
\ref{lem:7.4.4}.
Now assume that $\oplus_{\beta \in S^i} V^\beta$ has a structure of 
a simple VOA for $1\leq i\leq r-1$.
Then the set 
$$
  \mathcal{T}=\{\, V^\beta,\ V^{\beta+\alpha_{i+1}} \mid \beta\in S^i\,\}
$$
consists of inequivalent simple current $V^0$-modules with $(S^i\oplus \Z_2)
=S^{i+1}$-graded fusion rules:
$$
\begin{array}{l}
  V^{\beta_1}\times V^{\beta_2}
  = V^{\beta_1+\beta_2},\q 
  V^{\beta_1}\times V^{\beta_2+\alpha_{i+1}}
  = V^{\beta_1+\beta_2+\alpha_{i+1}},\q
  V^{\beta_1+\alpha_{i+1}}\times V^{\beta_2+\alpha_{i+1}}
  = V^{\beta_1+\beta_2},
\end{array}
$$
where $\beta_1,\beta_2\in S^i$.
By inductive assumption, $\oplus_{\beta\in S^i}V^\beta$ is an $S^i$-graded
simple current extension of $V^0$, and by Lemma \ref{lem:7.4.4}, a direct sum 
$V^0\oplus V^{\alpha_{i+1}}$ becomes a $\Z_2$-graded simple current extension 
of $V^0$.
Therefore, we can apply Theorem \ref{extension} to $\mathcal{T}$ 
to obtain the $S^{i+1}$-graded simple current extension 
$\oplus_{\beta\in S^{i+1}} V^\alpha$ of $V^0$.
Repeating this procedure, we finally obtain $S^r=S$-graded simple 
current extension $V=\oplus_{\alpha \in S} V^\alpha$ of $V^0=U_D$.
\qed

\begin{rem}
  In \cite{M4}, Miyamoto assumed stronger conditions than that in 
  Hypothesis \ref{I}. 
  In particular, he assumed that the structure codes $(D,S)$ are 
  of length $8k$ for some positive integer $k$. 
  Our refinement enable us to construct 2A-framed VOAs with structure
  codes of any length as long as Hypothesis \ref{I} is satisfied.
\end{rem}

\paragraph{Extension to SVOA.}
Let $V$ be a 2A-framed VOA with structure codes $(D,S)$.
Then by definition we have a decomposition $V=\oplus_{\alpha \in S}V^\alpha$ 
such that $V^0\simeq U_D$ and $\tau(V^\alpha)=\alpha$.
Assume that the pair $(D,S)$ satisfies the condition (1) of 
Hypothesis \ref{I}.
Then $V$ is an $S$-graded simple current extension of $V^0$ by Lemma 
\ref{lem:7.4.3}.
Suppose that there is a vector $\gamma \in S^\perp \setminus D$ such that 
$\la \gamma,\gamma \ra=1$.
Since the powers of $z$ in an $L(1/2,0)$-intertwining operator of type
$L(1/2,1/2)\times L(1/2,1/16)\to L(1/2,1/16)$ are half-integral, 
the powers of $z$ in a $U_D$-intertwining operator of type $U_{D+\gamma}
\times V^\alpha \to U_{D+\gamma}\fusion_{U_D} V^\alpha$ 
are integral for all $\alpha \in S$.
Therefore, all $V^\alpha \fusion_{U_D} U_{D+\gamma}$, $\alpha \in S$, 
have half-integral top weights.
Since the 1/16-word of $V^\alpha \fusion_{U_D} U_{D+\gamma}$ is $\alpha$ by 
the fusion rules \eqref{fusion ising}, the induced module 
$$
  \ind_{V^0}^{V} U_{D+\gamma}
  =\bigoplus_{\alpha \in S} V^\alpha \fusion_{V^0} U_{D+\gamma}
$$
has a unique $V$-module structure by Theorem \ref{induced module}.
Since $U_D\oplus U_{D+\gamma}$ is a simple current super-extension of 
$V^0=U_D$ by Theorem \ref{code VOA}, by applying Theorem \ref{extension}
we obtain:

\begin{thm}\label{super}
  With reference to the setup above, $V\oplus \ind_{V^0}^V U_{D+\gamma}$
  forms a simple current super-extension of $V$.
\end{thm}

\section{The baby-monster SVOA}

As shown in \cite{DMZ}, the moonshine VOA $V^\natural$ has a 
2A-frame $\w^\natural=e^1+\cds +e^{48}$.
One of its structure codes are determined in \cite{DGH} and \cite{M5}.
Let $S$ be the Reed-M\"{u}ller code $RM(4,1)$ defined as follows:
$$
  S=\Span_{\Z_2}\{ 
    (1^{16}),\ 
    (1^80^8),\ 
    (1^40^41^40^4),\ 
    (1^20^21^20^21^20^21^20^2),\ 
    (1010101010101010) 
  \} \subset \Z_2^{16}.
$$
Then define 
$$
  S^\natural:=\{ 
    (\alpha,\alpha,\alpha),\ 
    (\alpha^c,\alpha,\alpha),\ 
    (\alpha,\alpha^c,\alpha),\ 
    (\alpha,\alpha,\alpha^c) \in \Z_2^{48} 
    \mid \alpha \in S,\ \alpha^c:= \alpha +(1^{16})
  \} 
$$
and $D^\natural:=(S^\natural)^\perp$.

\begin{thm}(\cite{DGH} \cite{M5})
  The moonshine VOA $V^\natural$ has a structure codes 
  $(D^\natural,S^\natural)$.
\end{thm}

One can easily check that the pair $(D^\natural,S^\natural)$ satisfies 
the condition (1) of Hypothesis \ref{I}.
Thus, by Lemma \ref{lem:7.4.3}, we have

\begin{cor}
  Let $V^\natural=\oplus_{\alpha \in S^\natural} (V^\natural)^\alpha$
  be the 1/16-word decomposition according to the structure codes 
  $(D^\natural,S^\natural)$.
  Then the pair $(D^\natural,S^\natural)$ and the set 
  $\{ (V^\natural)^\alpha \mid \alpha \in S^\natural\}$ 
  satisfy Hypothesis \ref{I}.
  Therefore, $V^\natural$ is an $S^\natural$-graded simple current 
  extension of $(V^\natural)^0=U_{D^\natural}$.
\end{cor}

Now set $e=e^1$ and consider the commutant subalgebra $T_e^\natural(0)$
of $\vir (e)$ in $V^\natural$.
Since $\{ 1\} \cap \supp (S^\natural)\ne \emptyset$, 
$V^\natural_e(1/16)$ is not zero.
Then by the condition (1) of Hypothesis \ref{I}, $V_e(1/2)$ is not zero,
too.
Therefore, we obtain a decomposition
$$
  V^\natural
  = L(1/2,0)\tensor T_e^\natural(0) 
  \oplus L(1/2,1/2)\tensor T_e^\natural (1/2)
  \oplus L(1/2,1/16)\tensor T_e^\natural (1/16)
$$
such that $T_e^\natural(h)\ne 0$ for $h=0,1/2,1/16$.
By Theorem \ref{SVOA}, we know that $T_e^\natural(0)\oplus 
T_e^\natural(1/2)$ has a structure of a simple vertex operator 
superalgebra.
This algebra was first considered by H\"{o}hn \cite{H1} and he called it 
the {\it baby-monster SVOA}, because the centralizer 
$C_{\aut (V^\natural)}(\tau_e)$ is isomorphic to the 2-fold 
central extension $\la \tau_e\ra\cd \B$  of the baby-monster 
sporadic finite simple group $\B$ \cite{ATLAS} and so 
$\B$ naturally acts on it.
Following him, we set $\VB^0:=T_e^\natural(0)$, 
$\VB^1:=T_e^\natural(1/2)$ and $\VB:=T_e^\natural(0)\oplus 
T_e^\natural(1/2)$.
We also know that $T_e^\natural(1/16)$ is an irreducible $\Z_2$-twisted
$\VB$-module by Theorem \ref{twisted piece}, and so we set 
$\VB_T :=T_e^\natural(1/16)$ for convention.
Since all the conformal vectors of $V^\natural$ with central charge $1/2$
are conjugate under the Monster $\M =\aut (V^\natural)$ by \cite{C} and 
\cite{M1}, the algebraic structures on $\VB$ and $\VB_T$ are independent 
of choice a conformal vector $e=e^1$.

By definition, the Virasoro vector of $\VB^0$ is given by 
$\w^\natural -e^1=e^2+\cds +e^{48}$.
Thus $\VB^0$ is a 2A-framed VOA.
We compute the structure codes of $\VB^0$.
Set 
$$
  (S^\natural)^0:= \{ \alpha \in S^\natural \mid \{ 1\} \cap 
  \supp (\alpha) =\emptyset\},\q
  (S^\natural)^1:= \{ \alpha \in S^\natural \mid \{ 1\} \cap
  \supp (\alpha) =\{ 1\}\} .
$$
Then $S^\natural=(S^\natural)^0\sqcup (S^\natural)^1$ and by definition 
of Miyamoto involution, $(V^\natural)^{\la \tau_{e}\ra}
=\oplus_{\alpha \in (S^\natural)^0} (V^\natural)^\alpha$ and 
$V_e^\natural (1/16)= \oplus_{\beta \in (S^\natural)^1} (V^\natural)^\beta$.
Since $V^\natural$ is an $S^\natural$-graded simple current extension
of $(V^\natural)^0$, the fixed point subalgebra 
$(V^\natural)^{\la \tau_e\ra}$ is also an $(S^\natural)^0$-graded 
simple current extension of $(V^\natural)^0$ and $V_e^\natural(1/16)$ is 
a simple current $(V^\natural)^{\la \tau_e\ra}$-module by 
Corollary \ref{cor:2.13}.
Now define $\phi_\epsilon: \Z_2^{47}\hookrightarrow \Z_2^{48}$ by 
$\Z_2^{47} \ni \alpha\mapsto (\epsilon,\alpha)\in \Z_2^{48}$ for 
$\epsilon=0,1$, and set 
$$
  D^{\flat,\epsilon}:=\{ \alpha \in \Z_2^{47} \mid  \phi_\epsilon (\alpha) 
  \in D^\natural\} , \  \epsilon =0,1, 
  \q
  S^\flat:=\{ \beta \in \Z_2^{47} \mid \phi_0 (\beta) \in (S^\natural)^0\} .
$$

\begin{prop}
  The structure codes of $\VB^0$ with respect to the 2A-frame 
  $e^2+\cds +e^{48}$ are $(D^{\flat,0},S^\flat)$.
\end{prop}

\pf
For $\alpha \in (S^\natural)^0$, define $(V^\natural)^{\alpha,\epsilon}$
to be the sum of all irreducible $\tensor_{i=1}^{48} \vir(e^i)$-submodules 
of $(V^\natural)^\alpha$ whose $\vir (e^1)$-components are isomorphic to 
$L(1/2,\epsilon/2)$ for $\epsilon =0,1$.
Then $V_e^\natural(1/2)\ne 0$ implies that $(V^\natural)^{\alpha,\epsilon}
\ne 0$ for all $\alpha \in (S^\natural)^0$ and $\epsilon =0,1$.
Therefore, $(V^\natural)^\alpha = (V^\natural)^{\alpha,0}\oplus 
(V^\natural)^{\alpha,1}$ and we obtain 1/16-word decompositions 
$V_e^\natural(0)=\oplus_{\alpha \in (S^\natural)^0} 
(V^\natural)^{\alpha,0}$ 
and 
$V_e^\natural(1/2)=\oplus_{\alpha \in (S^\natural)^0} 
(V^\natural)^{\alpha,1}$.
Since $D^\natural=\phi_0(D^{\flat,0})\sqcup \phi_1(D^{\flat,1})$, 
$(V^\natural)^{0,0}\simeq L(1/2,0)\tensor U_{D^{\flat,0}}$.
Thus $\VB^0$ has a 1/16-word decomposition $\VB^0=\oplus_{\alpha 
\in S^\flat} (\VB^0)^\alpha$ such that $\tau((\VB^0)^\alpha)=
\alpha$ and $(\VB^0)^0\simeq U_{D^{\flat,0}}$.
Hence the structure codes of $\VB^0$ are $(D^{\flat,0},S^\flat)$.
\qed 

\begin{rem}\label{rem:5.4}
  By the proof above, we find that $\VB^1$ also has a 1/16-word decomposition 
  $\VB^1=\oplus_{\alpha \in S^\flat} (\VB^1)^\alpha$ 
  such that $\tau((\VB^1)^\alpha)=\alpha$.
  In particular, $(\VB^1)^0$ is isomorphic to a coset module 
  $U_{D^{\flat,1}}$.
\end{rem}

The following is easy to see:

\begin{lem}
  The pair $(D^{\flat,0},S^\flat)$ satisfies the condition (1) of 
  Hypothesis \ref{I}.
\end{lem}

Therefore, $\VB^0=\oplus_{\alpha \in S^\flat} (\VB^0)^\alpha$ is an 
$S^\flat$-graded simple current extension of the code VOA 
$U_{D^{\flat,0}}$.
Since the pair $(D^{\flat,0},S^\flat)$ and the set 
$\{ (\VB^0)^\alpha \mid \alpha \in S^\flat\}$ satisfy Hypothesis 
\ref{I}, we can construct $\VB^0$ without reference to $V^\natural$
by Theorem \ref{M4main}.

\begin{prop}
  $\VB^1$ is a simple current $\VB^0$-module.
\end{prop}

\pf
By Remark \ref{rem:5.4}, $\VB^1$ is isomorphic to the induced module 
$\ind_{U_{D^{\flat,0}}}^{\VB^0} U_{D^{\flat,1}}$ which is an
$S^\flat$-stable $\VB^0$-module.
Therefore, we have the fusion rule 
$$
  \VB^1\times \VB^1= \ind_{U_{D^{\flat,0}}}^{\VB^0} (U_{D^{\flat,1}}
  \fusion_{U_{D^{\flat,0}}} U_{D^{\flat,1}})
  = \ind_{U_{D^{\flat,0}}}^{\VB^0} U_{D^{\flat,0}}
  =\VB^0
$$
by Theorem \ref{lifting}.
Thus $\VB^1$ is a simple current $\VB^0$-module by Lemma \ref{criterion}.
\qed

\begin{rem}
  We note that by using Theorem \ref{super} we can define the SVOA structure 
  on $\VB$ without reference to $V^\natural$.
\end{rem}

Up to now, we have established that $V^\natural$ and its conformal 
vector $e$ satisfy all the conditions in Hypothesis \ref{II}.
Moreover, it is shown in \cite{C} and \cite{M1} that there is a one-to-one 
correspondence between the set of conformal vectors of $V^\natural$ with 
central charge 1/2 and the set of corresponding Miyamoto involutions
on $V^\natural$ which is known to be the 2A-conjugacy class of the Monster.
Therefore, $C_{\aut (V^\natural)}(e)=C_{\M}(\tau_e)\simeq \la \tau_e\ra
\cd \B$.
Thus, $C_{\aut (V^\natural)}(e)/\la \tau_e\ra$ is a simple group and 
we can apply Theorem \ref{one-point stabilizer} to $\VB$.

\begin{thm}\label{baby}
  (1)\ The SVOA $\VB$ obtained from $V^\natural$ by
  cutting off the Ising model is a simple SVOA.
  \\
  (2)\ The piece $\VB_T$ obtained from $V^\natural$ is an
  irreducible $\Z_2$-twisted $\VB$-module.
  \\
  (3)\ $\aut (\VB^0)\simeq \mathbb{B}$ and $\aut (\VB)\simeq 2\times 
  \mathbb{B}$.
  \\
  (4)\ $\VB_T$ as a $\VB^0$-module is irreducible. 
  Thus, there are exactly three irreducible $\VB^0$-modules,
  $\VB^0$, $\VB^1$ and $\VB_T$.
  \\
  (5)\ The fusion rules for irreducible $\VB^0$-modules are as follows:
  $$
    \VB^1\times \VB^1=\VB^0,\q \VB^1\times \VB_T=\VB_T,\q 
    \VB_T\times \VB_T=\VB^0+\VB^1.
  $$
\end{thm}

\pf
(1) follows from Theorem \ref{SVOA}, (2) follows from 
\ref{twisted piece} and (3) will follow from Theorem 
\ref{one-point stabilizer} and the fact $C_{\aut (V^\natural)}(e)
=\la \tau_e\ra\cd \B$.

Consider (4).
By Corollary \ref{cor:3.15}, $\VB_T$ as a $\VB^0$-module is irreducible.
Then the assertion follows from Proposition \ref{prop:8.1.9}.

Consider (5).
We only have to show the fusion rule $\VB_T \times \VB_T=\VB^0+\VB^1$.
By considering the 1/16-word of $\VB_T$, the fusion product 
$\VB_T\times \VB_T$ is contained in $\N \VB^0\oplus \N \VB^1$ in the fusion 
algebra for $\VB^0$.
Write $\VB_T\times \VB_T=n_0 \VB^0+n_1\VB^1$ with $n_0,n_1 \in \N$.
Then the simplicity of $V^\natural$ implies $n_0\ne 0$ and $n_1\ne 0$.
And by applying $\VB^1$ to $\VB_T\times \VB_T$, we see that $n_0=n_1$.
Since the dual module of $\VB_T$ is isomorphic to $\VB_T$, 
the space of $\VB^0$-intertwining operator of type $\VB_T\times 
\VB_T\to \VB^0$ is one-dimensional.
Thus $n_0=n_1=1$ as desired.
\qed

\begin{rem}
  The assertion (1) of Theorem \ref{baby} is already shown by H\"{o}hn 
  in \cite{Ho1}, and (3) of Theorem \ref{baby} is also proved in 
  \cite{Ho2}.
  However, H\"{o}hn's proofs in \cite{Ho1} \cite{Ho2} and ours are 
  quite different.
  In particular, in \cite{Ho2}, he used many results on the baby-monster 
  simple group.
  Our argument can be applied to any 2A-framed VOAs satisfying Hypothesis 
  \ref{II} and Hypothesis \ref{I} since we have only used the facts that 
  $C_{\aut (V^\natural)}(e)= C_{\aut (V^{\natural})}(\tau_e)$ and 
  $C_{\aut (V^\natural)}(\tau_e)/\la \tau_e\ra$ is a simple group. 
\end{rem}

The classification of irreducible $\VB^0$-modules has many interesting
corollaries.

\begin{cor}
  The irreducible 2A-twisted $V^\natural$-module has a shape
  $$
    L(1/2,1/2)\tensor \VB^0
    \oplus L(1/2,0)\tensor \VB^1
    \oplus L(1/2,1/16)\tensor \VB_T.
  $$
\end{cor}

\pf
Follows from Theorem \ref{baby}, Theorem \ref{thm:8.1.8} and 
Proposition \ref{prop:8.1.9}.
\qed

\begin{rem}
  A straightforward construction of the 2A-twisted (and 2B-twisted) 
  $V^\natural$-module is given by Lam \cite{L2}.
  In his construction, it is given as $U_{D^\natural+\gamma}
  \fusion_{U_{D^\natural}} V^\natural$ with $\gamma=(10^{47})\in 
  (\Z/2\Z)^{48}$.
\end{rem}

\begin{cor}
  For any conformal vector $e\in V^\natural$ with central charge 1/2,
  there is no automorphism $\rho$ on $V^\natural$ such that
  $\rho (V_e^\natural(h))=V_e^\natural(h)$ for $h=0,1/2$ and 
  $\rho|_{(V^\natural)^{\la \tau_e\ra}}=\sigma_e$.
\end{cor}

\pf
Suppose such an automorphism $\rho$ exists.
We remark that $\rho$ also preserves the space $V_e(1/16)$ as $\rho \in 
C_{\aut (V)}(e)$.
We view $V_e^\natural (1/16)$ as a $(V^\natural)^{\la \tau_e\ra}$-module 
by a restriction of the vertex operator map $Y_{V^\natural}(\cd,z)$ 
on $V^\natural$.
Consider the $\sigma_e$-conjugate $(V^\natural)^{\la \tau_e\ra}$-module 
$V^\natural_e(1/16)^{\sigma_e}$.
By Theorem \ref{baby} and Proposition \ref{prop:8.1.9}, 
$V_e^\natural(1/16)^{\sigma_e}$ is not isomorphic to 
$V_e^\natural(1/16)$ as a $(V^\natural)^{\la \tau_e\ra}$-module.
On the other hand, we can take a canonical linear isomorphism  
$\varphi: V_e^\natural (1/16)\to V_e^\natural (1/16)^{\sigma_e}$ 
such that $Y_{V^\natural_e(1/16)^{\sigma_e}}(a,z) \varphi v
=\varphi Y_{V^\natural}(\sigma_e a,z) v$ for all $a\in 
(V^\natural)^{\la \tau_e\ra}$ and $v\in V_e^\natural(1/16)$ 
by definition of the conjugate module.
Then we have 
$$ 
  Y_{V_e^\natural(1/16)^{\sigma_e}}(a,z)\varphi \rho v
  = \varphi Y_{V^\natural}(\sigma_e a,z) \rho v
  = \varphi Y_{V^\natural}(\rho a,z) \rho v
  = \varphi \rho Y_{V^\natural}(a,z)v 
$$
for any $a\in (V^\natural)^{\la \tau_e\ra}$ and $v\in V_e^\natural(1/16)$.
Thus $\varphi \rho$ defines a $(V^\natural)^{\la \tau_e\ra}$-isomorphism
between $V_e^\natural(1/16)$ and $V_e^\natural(1/16)^{\sigma_e}$, 
which is a contradiction.
\qed

\begin{cor}
  The 2A-orbifold construction applied to the moonshine VOA $V^\natural$
  yields $V^\natural$ itself again.
\end{cor}

\pf
Follows from Theorem \ref{baby} and Corollary \ref{cor:3.15}.
\qed

\begin{rem}
  The statement in the corollary above was conjectured by Tuite 
  \cite{Tu}.
  In \cite{Tu}, Tuite has shown that any $\Z_p$-orbifold construction
  of $V^\natural$ yields the moonshine VOA $V^\natural$ or the Leech 
  lattice VOA $V_\Lambda$ under the uniqueness conjecture of 
  the moonshine VOA which states that $V^\natural$ constructed by Frenkel
  et.\ al.\ \cite{FLM} is the unique holomorphic VOA with central charge 24 
  whose weight one subspace is trivial.
\end{rem}

Finally, we end this paper by presenting the modular transformations 
of characters of $\VB^0$-modules.
Here the character means the conformal character, not the $q$-dimension, 
of modules.
Recall the characters of $L(1/2,0)$-modules.
By our explicit construction in Section \ref{Ising}, 
one can easily prove the following (cf.\ \cite{FFR} \cite{FRW}):
$$
\begin{array}{l}
  \ch_{L(1/2,0)}(\tau)
  = (1/2)\cd q^{-1/48}\l\{ \prod_{n=0}^\infty (1+q^{n+1/2})
    + \prod_{n=0}^\infty (1-q^{n+1/2})\r\} , 
  \vsv\\
  \ch_{L(1/2,1/2)}(\tau)
  = (1/2)\cd q^{-1/48}\l\{ \prod_{n=0}^\infty (1+q^{n+1/2})
    - \prod_{n=0}^\infty (1-q^{n+1/2})\r\} , 
  \vsv\\
  \ch_{L(1/2,1/16)}(\tau)= q^{-1/24} \prod_{n=1}^\infty (1+q^n).
\end{array}
$$
The following modular transformations are well-known:
$$
\begin{array}{lll}
  \ch_{L(1/2,0)}(-1/\tau) 
    &=& \dfr{1}{2} \ch_{L(1/2,0)}(\tau)
    +\dfr{1}{2} \ch_{L(1/2,1/2)}(\tau)
    +\dfr{1}{\sqrt{2}} \ch_{L(1/2,1/16)}(\tau), 
  \vsb\\
  \ch_{L(1/2,1/2)}(-1/\tau) 
    &=& \dfr{1}{2} \ch_{L(1/2,0)}(\tau)
    +\dfr{1}{2} \ch_{L(1/2,1/2)}(\tau)
    -\dfr{1}{\sqrt{2}} \ch_{L(1/2,1/16)}(\tau),
  \vsb\\
  \ch_{L(1/2,1/16)}(-1/\tau)
    &=& \dfr{1}{\sqrt{2}} \ch_{L(1/2,0)}(\tau)
    -\dfr{1}{\sqrt{2}} \ch_{L(1/2,1/2)}(\tau).
\end{array}
$$
Set $j(\tau):=J(\tau)-744$, where $J(\tau)$ is the famous 
$\SL_2(\Z)$-invariant.
Since $\ch_{V^\natural}(\tau)=j(\tau)$ and 
$$
  \ch_{V^\natural}(\tau)
  =\ch_{L(1/2,0)}(\tau) \ch_{\VB^0}(\tau)
  +\ch_{L(1/2,1/2)}(\tau) \ch_{\VB^1}(\tau)
  +\ch_{L(1/2,1/16)}(\tau) \ch_{\VB_T}(\tau),
$$
we can write down the characters of irreducible $\VB^0$-modules 
by using those of $V^\natural$ and $L(1/2,0)$-modules.
This computation is already done in \cite{Ma} by using 
Norton's trace formula.
The results are written as a rational expression involving the functions
$j(\tau)$, $\ch_{L(1/2,h)}(\tau)$, $h=0,1/2,1/16$, their first and 
second derivatives and the Eisenstein series $E_2(\tau)$ and 
$E_4(\tau)$, see \cite{Ma}.

By Zhu's theorem \cite{Z}, the linear space spanned by 
$\{ \ch_{\VB^0}(\tau),\ch_{\VB^1}(\tau),\ch_{\VB_T}(\tau)\}$ 
affords an $\SL_2(\Z)$-action.
By using the modular transformations
for $j(\tau)$ and $\ch_{L(1/2,h)}(\tau)$, $h=0,1/2,1/16$, we can 
show the following modular transformations:
$$
\begin{array}{lll}
  \ch_{\VB^0}(-1/\tau) 
    &=& \dfr{1}{2} \ch_{\VB^0}(\tau)
    +\dfr{1}{2}    \ch_{\VB^1}(\tau)
    +\dfr{1}{\sqrt{2}} \ch_{\VB_T}(\tau), 
  \vsb\\
  \ch_{\VB^1}(-1/\tau) 
    &=& \dfr{1}{2} \ch_{\VB^0}(\tau)
    +\dfr{1}{2} \ch_{\VB^1}(\tau)
    -\dfr{1}{\sqrt{2}} \ch_{\VB_T}(\tau),
  \vsb\\
  \ch_{\VB_T}(-1/\tau)
    &=& \dfr{1}{\sqrt{2}} \ch_{\VB^0}(\tau)
    -\dfr{1}{\sqrt{2}} \ch_{\VB^1}(\tau).
\end{array}
$$
Namely, we have exactly the same modular transformation laws
for the Ising model $L(1/2,0)$.
As in Theorem \ref{baby}, we also note that the fusion algebra 
for $\VB^0$ is also canonically isomorphic to that of $L(1/2,0)$.
Therefore, we may say that $L(1/2,0)$ and $\VB^0$ form a dual-pair
in the moonshine VOA $V^\natural$.


\begin{quote}
\begin{center}
  {\bf Acknowledgments}
\end{center}
  The author wishes to thank the members of the Komaba Seminar on 
  Finite Groups for valuable discussions and for noticing a gap in 
  a prototype version of this paper.
  He also thank Professor Atsushi Matsuo for the information of Tuite's 
  work and Professor Ching Hung Lam for his suggestions improving 
  this paper.
\end{quote}

\end{document}